\pdfoutput=1
\RequirePackage{ifpdf}
\ifpdf 
\documentclass[pdftex]{sigma}
\else
\documentclass{sigma}
\fi

\usepackage{tikz}
\usetikzlibrary{matrix}
\usepackage{mathtools}
\newtheorem{Theorem}{Theorem}[section]
\newtheorem{Corollary}[Theorem]{Corollary}
\newtheorem{Proposition}[Theorem]{Proposition}
{ \theoremstyle{definition}
\newtheorem{Definition}[Theorem]{Definition}
\newtheorem{Remark}[Theorem]{Remark} }

\numberwithin{equation}{section}

\begin{document}

\newcommand{\arXivNumber}{1312.0909}

\allowdisplaybreaks

\renewcommand{\PaperNumber}{071}

\FirstPageHeading

\ShortArticleName{Spherical Functions of Fundamental $K$-Types Associated with the~$n$-Dimensional Sphere}

\ArticleName{Spherical Functions of Fundamental $\boldsymbol{K}$-Types\\
Associated with the $\boldsymbol{n}$-Dimensional Sphere}

\Author{Juan Alfredo TIRAO and Ignacio Nahuel ZURRI\'AN}

\AuthorNameForHeading{J.A.~Tirao and I.N.~Zurri\'an}

\Address{CIEM-FaMAF, Universidad Nacional de C\'or\-do\-ba, Argentina}
\Email{\href{mailto:tirao@famaf.unc.edu.ar}{tirao@famaf.unc.edu.ar}, \href{mailto:zurrian@famaf.unc.edu.ar}{zurrian@famaf.unc.edu.ar}}
\URLaddress{\url{http://www.famaf.unc.edu.ar/~zurrian/}}

\ArticleDates{Received December 20, 2013, in f\/inal form June 20, 2014; Published online July 07, 2014}

\Abstract{In this paper, we describe the irreducible spherical functions of fundamental $K$-types associated with the
pair $(G,K)=({\mathrm{SO}}(n+1),{\mathrm{SO}}(n))$ in terms of matrix hyper\-geo\-metric functions.
The output of this description is that the irreducible spherical functions of the same $K$-fundamental type are encoded
in new examples of classical sequences of matrix-valued orthogonal polynomials, of size~$2$ and~$3$, with respect to
a~matrix-weight~$W$ supported on~$[0,1]$.
Moreover, we show that~$W$ has a~second order symmetric hypergeometric operator~$D$.}


\Classification{22E45; 33C45; 33C47}

\section{Introduction}

The theory of spherical functions dates back to the classical papers of \'E.~Cartan and H.~Weyl; they showed that spherical harmonics arise in a~natural way from the study of functions on
the~$n$-dimensional sphere $S^n={\mathrm{SO}}(n+1)/{\mathrm{SO}}(n)$.
The f\/irst general results in this direction were obtained in 1950 by Gel'fand, who considered zonal spherical functions
of a~Riemannian symmetric space $G/K$.
In this case we have a~decomposition $G=KAK$.
When the Abelian subgroup~$A$ is one dimensional, the restrictions of zonal spherical functions to~$A$ can be identif\/ied
with hypergeometric functions, providing a~deep and fruitful connection between group representation theory and special
functions.
In particular when~$G$ is compact this gives a~one to one correspondence between all zonal spherical functions of the
symmetric pair $(G,K)$ and a~sequence of orthogonal polynomials.

In light of this remarkable background it is reasonable to look for an extension of the above results, by considering
matrix-valued irreducible spherical functions on~$G$ of a~general $K$-type.
This was accomplished for the f\/irst time in the case of the complex projective plane
$P_2(\mathbb{C})=\mathrm{SU}(3)/\mathrm{U}(2)$ in~\cite{GPT02a}.
This seminal work gave rise to a~series of papers including~\cite{GPT02b, GPT03, GPT05,HP13,KPR12,KPR13, PR08, PT04, PT07b,
PTZ13}, where one considers matrix valued spherical functions associated to a~compact symmetric pair
$(G,K)$ of rank one, arriving at sequences of matrix valued orthogonal polynomials of one real variable satisfying an
explicit three-term recursion relation, which are also eigenfunctions of a~second order matrix dif\/ferential operator
(bispectral property).

The very explicit results contained in this paper are obtained for certain $K$-types, namely the fundamental $K$-types.
Also, the detailed construction of sequences of matrix orthogonal polynomials out of these irreducible spherical
functions, following the general pattern established in~\cite{GPT02a}, gives new examples of {\it classical sequences}
of matrix-valued orthogonal polynomials of size~$2$ and~$3$.
For the general notions concerning matrix-valued orthogonal polynomials see~\cite{GT07}.
Interesting generalizations of these sequences are given in~\cite{PZ13}, where the coef\/f\/icients of the three term
recursion relation satisf\/ied by them is exhibited.

The present paper is an outgrowth of the results of~\cite[Chapter~5]{Z13} and we are currently working on the
extension of these results for the spherical functions of any $K$-type associated with the~$n$-dimensional sphere.
Using~\cite{TZ13a}, one can obtain the corresponding results for the spherical functions of any $K$-type associated
with~$n$-dimensional real projective space.
The starting point is to describe the irreducible spherical functions associated with the pair
$(G,K)=({\mathrm{SO}}(n+1),{\mathrm{SO}}(n))$ in terms of eigenfunctions of a~matrix linear dif\/ferential operator of
order two.
The output of this description is that the irreducible spherical functions of the same fundamental $K$-type are encoded
in a~sequence of matrix valued orthogonal polynomials.

Brief\/ly the main results of this paper are the following.
After some preliminaries, in Section~\ref{sec: delta} we study the eigenfunctions of an operator~$\Delta$ on~$G$, which
is closely related to the Casimir operator.
Every spherical function~$\Phi$ has to be eigenfunction of this operator~$\Delta$; considering the $KAK$-decomposition
\begin{gather*}
{\mathrm{SO}}(n+1)={\mathrm{SO}}(n){\mathrm{SO}}(2){\mathrm{SO}}(n)
\end{gather*}
and choosing an appropriate coordinate~$y$ on an open subset of~$A$, we translate the condition
$\Delta\Phi=\lambda\Phi$, $\lambda\in\mathbb{C}$, into a~matrix valued dif\/ferential equation $\widetilde DH=\lambda H$
on the open interval $(0,1)$, where~$H$ is the restriction of~$\Phi$ to ${\mathrm{SO}}(2)$.
The property of the spherical functions
\begin{gather*}
\Phi(xgy)=\pi(x)\Phi(g)\pi(y),
\qquad
g\in G,
\qquad
x,y\in K,
\end{gather*}
tell us that~$\Phi$ is determined by its $K$-type and the function~$H$.

In Section~\ref{mirred} we f\/irst explicitly describe all the irreducible spherical functions of the symmetric pair
$(G,K)=({\mathrm{SO}}(n+1),{\mathrm{SO}}(n))$ with~$M$-irreducible $K$-types, with $M={\mathrm{SO}}(n-1)$, the
centralizer of the subgroup~$A$ in~$K$; we give these expressions in terms of the hypergeometric function $_2F_1$.

In Section~\ref{sec: fund} the operator $\widetilde D$ is studied in detail when the $K$-types correspond to fundamental
representations.
Certain $K$-fundamental types are~$M$ irreducible, and therefore they were already considered en Section~\ref{mirred};
besides, when~$n$ is odd there is a~particular fundamental $K$-type which has three~$M$-submodules, this case is studied
in the last section of this work.
For the rest of the cases we considered separately when~$n$ is even and when~$n$ is odd.
Although, in both cases we worked with the concrete realizations of the fundamental representations considering the
exterior powers of the standard representation of ${\mathrm{SO}}(n)$:
\begin{gather*}
\Lambda^1\big(\mathbb{C}^n\big),
\
\Lambda^2\big(\mathbb{C}^n\big),
\
\dots,
\
\Lambda^{\ell-1}\big(\mathbb{C}^n\big),
\end{gather*}
with $n=2\ell$ or $n=2\ell+1$.

In Section~\ref{s2l} we conjugate the operator $\widetilde D$, by using the polynomial function
\begin{gather*}
\Psi(y)=\left(
\begin{matrix}
2y-1&1
\\
1&2y-1
\end{matrix}
\right),
\end{gather*}
whose columns correspond to irreducible spherical functions, in order to obtain a~matrix-valued hypergeometric operator
$D=\Psi^{-1}\widetilde D \Psi$:
\begin{gather*}
DP=y(1-y)P''+(C-yU)P'-VP,
\end{gather*}
with
\begin{gather*}
C=
\begin{pmatrix}
(n/2+1)&1
\\
1&(n/2+1)
\end{pmatrix},
\qquad
U=(n+2)I,
\qquad
V= \left(
\begin{matrix}
p&0
\\
0&n-p
\end{matrix}
\right).
\end{gather*}

Then, we study all the possible eigenvalues corresponding to irreducible spherical functions and all the polynomial
eigenfunctions of~$D$.

In Section~\ref{innerprod}, for any fundamental $K$-type ($\Lambda^k(\mathbb{C}^n)$) with $1\le p\le \ell-1$, we f\/ind
a~matrix-weight~$W$, which is a~scalar multiple of
\begin{gather*}
W= (y(1- y))^{n/2-1} \left(
\begin{matrix}
p(2y-1)^2+n-p&n(2y-1)
\\
n(2y-1)&(n-p)(2y-1)^2+p
\end{matrix}
\right),
\end{gather*}
such that~$D$ is a~symmetric operator with respect to the inner product def\/ined among continuous vector-valued functions
on $[0,1]$~by
\begin{gather*}
\langle P_1,P_2 \rangle_W =\int_0^1 P_2^*(y)W(y)P_1(y)dy.
\end{gather*}
Also we prove that every spherical function gives a~vector polynomial eigenfunction~$P$ of~$D$.
Therefore we obtain the following explicit expression of~$P$ in terms of the matrix hypergeometric function for any
irreducible spherical function
\begin{gather*}
P(y)=\sum\limits_{j=0}^{w}\frac{y^j}{j!} [C;U;V+\lambda]_j P(0),
\end{gather*}
see Theorem~\ref{columns}.

In Section~\ref{mvop} for each pair $(n,p)$ we construct a~sequence of matrix orthogonal polynomials $\{P_w\}_{w\ge0}$
of size $2$ with respect to the weight function~$W$, which are eigenfunctions of the symmetric dif\/ferential
operator~$D$.
Namely,
\begin{gather*}
D P_w = P_w \left(\begin{matrix}
\lambda(w,0) & 0
\\
0 & \lambda(w,1)
\end{matrix}
\right),
\end{gather*}
where
\begin{gather*}
\lambda(w,\delta)=
\begin{cases}
-w(w+n+1)-p
&\text{if}  \quad
\delta=0,
\\
-w(w+n+1)-n+p
&\text{if}  \quad
\delta=1.
\end{cases}
\end{gather*}

Finally, in Section~\ref{PC} we develop the same techniques in order to obtain analogous results for irreducible
spherical functions of the particular $K$-fundamental type $\Lambda^\ell(\mathbb{C}^n)$ for which we have
three~$M$-submodules instead of only two.
This only occurs when~$n$ is of the form $2\ell+1$.

It is worth to notice that, unlike the other cases, the $3\times3$ matrix-weight built here does reduce to a~smaller
size.

\section{Preliminaries}

\subsection{Spherical functions}

Let~$G$ be a~locally compact unimodular group and let~$K$ be a~compact subgroup of~$G$.
Let $\hat K$ denote the set of all equivalence classes of complex f\/inite dimensional irreducible representations of~$K$;
for each $\delta\in \hat K$, let $\xi_\delta$ denote the character of~$\delta$, $d(\delta)$ the degree of~$\delta$,
i.e.~the dimension of any representation in the class~$\delta$, and $\chi_\delta=d(\delta)\xi_\delta$.
We shall choose once and for all the Haar measure $dk$ on~$K$ normalized by $\int_K dk=1$.

We shall denote by~$V$ a~f\/inite dimensional vector space over the f\/ield $\mathbb{C}$ of complex numbers and by of all
linear transformations of~$V$ into~$V$.
Whenever we refer to a~topology on such a~vector space we shall be talking about the unique Hausdorf\/f linear topology on
it.

\begin{Definition}
A~spherical function~$\Phi$ on~$G$ of type $\delta\in \hat K$ is a~continuous function on~$G$ with values in
${\operatorname{End}}(V)$ such that
\begin{enumerate}\itemsep=0pt
\item[i)] $\Phi(e)=I$ ($I$ is the identity transformation);

\item[ii)] $\Phi(x)\Phi(y)=\int_K \chi_{\delta}(k^{-1})\Phi(xky) dk$ for all $x,y\in G$.
\end{enumerate}
\end{Definition}

The reader can f\/ind a~number of general results in~\cite{T77} and~\cite{GV88}.
For our purpose it is appropriate to recall the following facts.

\begin{Proposition}
[\protect{\cite[Proposition 1.2]{T77}}]
\label{propesf}
If $\Phi:G\longrightarrow {\operatorname{End}}(V)$ is a~spherical function of type~$\delta$ then:
\begin{enumerate}\itemsep=0pt
\item[$i)$] $\Phi(k_1gk_2)=\Phi(k_1)\Phi(g)\Phi(k_2)$, for all $k_1,k_2\in K$, $g\in G$;
\item[$ii)$] $k\mapsto \Phi(k)$ is a~representation of~$K$ such that any irreducible subrepresentation belongs to~$\delta$.
\end{enumerate}
\end{Proposition}

Concerning the def\/inition, let us point out that the spherical function~$\Phi$ determines its type univocally
(Proposition~\ref{propesf}) and let us say that the number of times that~$\delta$ occurs in the representation $k\mapsto
\Phi(k)$ is called the height of~$\Phi$.

A spherical function $\Phi: G \longrightarrow {\operatorname{End}}(V)$ is called irreducible if~$V$ has no proper
subspace invariant by $\Phi(g)$ for all $g\in G$.

If~$G$ is a~connected Lie group, it is not dif\/f\/icult to prove that any spherical function $\Phi:G\longrightarrow
{\operatorname{End}}(V)$ is dif\/ferentiable ($C^\infty$), and moreover that it is analytic.
Let $D(G)$ denote the algebra of all left invariant dif\/ferential operators on~$G$ and let $D(G)^K$ denote the subalgebra
of all operators in $D(G)$ which are invariant under all right translations by elements in~$K$.

In the following proposition $(V,\pi)$ will be a~f\/inite dimensional representation of~$K$ such that any irreducible
subrepresentation belongs to the same class $\delta\in\hat K$.

\begin{Proposition}
A~function $\Phi:G\longrightarrow {\operatorname{End}}(V)$ is a~spherical function of type~$\delta$ if and only if
\begin{enumerate}\itemsep=0pt
\item[$i)$] $\Phi$ is analytic;
\item[$ii)$] $\Phi(k_1gk_2)=\pi(k_1)\Phi(g)\pi(k_2)$, for all $k_1,k_2\in K$, $g\in G$, and $\Phi(e)=I$;
\item[$iii)$] $[D\Phi](g)=\Phi(g)[D\Phi](e)$, for all $D\in D(G)^K$, $g\in G$.
\end{enumerate}
\end{Proposition}

Moreover, we have that the eigenvalues $[D\Phi](e)$, $D\in D(G)^K$, characterize the spherical functions~$\Phi$ as
stated in the following proposition.

\begin{Proposition}[\protect{\cite[Remark 4.7]{T77}}]
Let $\Phi,\Psi:G\longrightarrow {\operatorname{End}}(V)$ be two spherical functions on a~connected Lie group~$G$ of the
same type $\delta\in K$.
Then $\Phi=\Psi$ if and only if $(D\Phi)(e)=(D\Psi)(e)$ for all $D\in D(G)^K$.
\end{Proposition}

Let us observe that if $\Phi:G\longrightarrow {\operatorname{End}}(V)$ is a~spherical function, then $\Phi:D\mapsto
[D\Phi](e)$ maps $D(G)^K$ into ${\operatorname{End}}_K(V)$ (${\operatorname{End}}_K(V)$ denotes the space of all linear
maps of~$V$ into~$V$ which commutes with $\pi(k)$ for all $k\in K$) def\/ining a~f\/inite dimensional representation of the
associative algebra $D(G)^K$.
Moreover, the spherical function is irreducible if and only if the representation $\Phi: D(G)^K\longrightarrow
{\operatorname{End}}_K (V)$ is irreducible.
We quote the following result from~\cite{PTZ13}.
\begin{Proposition}
[\protect{\cite[Proposition~2.5]{PTZ13}}]
Let~$G$ be a~connected reductive linear Lie group.
Then the following properties are equivalent:
\begin{enumerate}\itemsep=0pt
\item[$i)$] $D(G)^K$ is commutative;
\item[$ii)$] every irreducible spherical function of $(G,K)$ is of height one.
\end{enumerate}
\end{Proposition}

In this paper the pair $(G,K)$ is $({\mathrm{SO}}(n+1),{\mathrm{SO}}(n))$.
Then, it is known that $D(G)^K$ is an Abelian algebra; moreover, $D(G)^K$ is isomorphic to $D(G)^G\otimes D(K)^K$ (see
in~\cite[Theorem~10.1]{K89} or~\cite{Co}),
where $D(G)^G$ (resp.\ $D(K)^K$) denotes the subalgebra of all operators in $D(G)$ (resp.~$D(K)$)
which are invariant under all right translations by elements in~$G$ (resp.~$K$).

An immediate consequence of this is that all irreducible spherical functions of our pair $(G,K)$ are of height one.

Spherical functions of type~$\delta$ (see in~\cite[Section~3]{T77}) arise in a~natural way upon considering
representations of~$G$.
If $g\mapsto U(g)$ is a~continuous representation of~$G$, say on a~f\/inite dimensional vector space~$E$, then
\begin{gather*}
P_\delta=\int_K \chi_\delta\big(k^{-1}\big)U(k) dk
\end{gather*}
is a~projection of~$E$ onto $P_\delta E=E(\delta)$.
If $P_\delta\ne0$ the function $\Phi:G\longrightarrow {\operatorname{End}}(E(\delta))$ def\/ined~by
\begin{gather}
\label{proj}
\Phi(g)a=P_\delta U(g)a,
\qquad
g\in G,
\qquad
a\in E(\delta),
\end{gather}
is a~spherical function of type~$\delta$.
In fact, if $a\in E(\delta)$ we have
\begin{gather*}
\Phi(x)\Phi(y)a = P_\delta U(x)P_\delta U(y)a=\int_K \chi_\delta\big(k^{-1}\big) P_\delta U(x)U(k)U(y)a dk
\\
\phantom{\Phi(x)\Phi(y)a}{}
 =\left(\int_K\chi_\delta\big(k^{-1}\big)\Phi(xky) dk\right) a.
\end{gather*}

If the representation $g\mapsto U(g)$ is irreducible then the associated spherical function~$\Phi$ is also irreducible.
Conversely, any irreducible spherical function on a~compact group~$G$ arises in this way from a~f\/inite dimensional
irreducible representation of~$G$.

\subsection[Root space structure of $\mathfrak{so}(n,\mathbb{C})$]{Root space structure of $\boldsymbol{\mathfrak{so}(n,\mathbb{C})}$}

Let $E_{ik}$ denote the square matrix with a~1 in the $ik$-entry and zeros elsewhere; and let us consider the matrices
\begin{gather*}
I_{ki}=E_{ik}-E_{ki},
\qquad
1\le i,k\le n.
\end{gather*}
Then, the set $\{I_{ki}\}_{i<k}$ is a~basis of the Lie algebra $\mathfrak{so}(n)$.
These matrices satisfy the following commutation relations
\begin{gather*}
[I_{ki},I_{rs}]=\delta_{ks}I_{ri}+\delta_{ri}I_{sk}+\delta_{is}I_{kr}+\delta_{rk}I_{is}.
\end{gather*}
If we assume that $k>i$, $r>s$ then we have
\begin{gather*}
[I_{ki},I_{is}]=I_{sk},
\qquad
[I_{ki},I_{rk}]=I_{ri},
\qquad
[I_{ki},I_{ri}]=I_{kr},
\qquad
[I_{ki},I_{ks}]=I_{is},
\end{gather*}
and all the other brackets are zero.
From this it easily follows that the set
\begin{gather*}
\{I_{p,p-1}:2\le p\le n\}
\end{gather*}
generates the Lie algebra $\mathfrak{so}(n)$.

\begin{Proposition}
\label{rightinv}
Given $n\in\mathbb{N}$, we have that the operator
\begin{gather*}
Q_n=\sum\limits_{1\le i,k\le n}I_{ki}^2\in D({\mathrm{SO}}(n))
\end{gather*}
is right invariant under ${\mathrm{SO}}(n)$, i.e.~
\begin{gather*}
Q_n\in D({\mathrm{SO}}(n))^{{\mathrm{SO}}(n)},
\qquad
\forall\, n\in\mathbb{N}_0.
\end{gather*}
\end{Proposition}
\begin{proof}
To prove that $Q_n$ is right invariant under ${\mathrm{SO}}(n)$ it is enough to prove that $\dot I_{p,p-1}(Q_n)=0$ for
all $2\le p\le n$.
We have
\begin{gather*}
\dot I_{p,p-1}(Q_n)=\sum\limits_{1\le i,k\le n}\big([I_{p,p-1},I_{ki}]I_{ki}+I_{ki}[I_{p,p-1},I_{ki}]\big).
\end{gather*}
Then
\begin{gather*}
\dot I_{p,p-1}(Q_n)= \sum\limits_{1\le i\le n}(I_{ip}I_{p-1,i}+I_{p-1,i}I_{ip})+\sum\limits_{1\le k\le n}(I_{k,p-1}I_{kp}+I_{kp}I_{k,p-1})
\\
\phantom{\dot I_{p,p-1}(Q_n)=}
{}+\sum\limits_{1\le k\le n}(I_{pk}I_{k,p-1}+I_{k,p-1}I_{pk})+\sum\limits_{1\le i\le n}(I_{p-1,i}I_{p,i}+I_{p,i}I_{p-1,i})=0.
\end{gather*}
This proves the proposition.
\end{proof}

\subsection[The operator $Q_{2\ell}$]{The operator $\boldsymbol{Q_{2\ell}}$}

Let us assume that $n=2\ell$.
We look at a~root space decomposition of $\mathfrak{so}(n)$ in terms of the basis elements $I_{ki}$, $1\le i<k\le n$.

The linear span
\begin{gather*}
\mathfrak h={\langle I_{21},I_{43},\dots,I_{2\ell,2\ell-1}\rangle}_\mathbb{C}
\end{gather*}
is a~Cartan subalgebra of $\mathfrak{so}(n,\mathbb{C})$.
To f\/ind the root vectors it is convenient to visualize the elements of $\mathfrak{so}(n,\mathbb{C})$ as $\ell\times\ell$
matrices of $2\times2$ blocks.
Thus $\mathfrak h$ is the subspace of all diagonal matrices of $2\times2$ skew-symmetric blocks.
The subspaces of all matrices~$A$ with a~block $A_{jk}$ of size two, $1\le j<k\le\ell$, in the place $(j,k)$ and
$-A_{jk}^t$ in the place $(k,j)$ with zeros in all other places, are $\operatorname{ad}(\mathfrak h)$-stable.
Let
\begin{gather*}
H=i(x_1I_{21}+\dots+x_\ell I_{2\ell,2\ell-1})\in\mathfrak h,
\end{gather*}
for $x_1,\dots,x_\ell\in\mathbb{R}$.
Then $[H,A]=\lambda(H)A$, $\forall\, H\in\mathfrak{h}$, if and only if for every $A_{jk}$ we have
\begin{gather*}
x_j(H)iI_{2j,2j-1}A_{jk}-x_k(H)iA_{jk}I_{2k,2k-1}=\lambda(H) A_{jk},
\qquad
\forall\, H\in \mathfrak{h}.
\end{gather*}
Up to a~scalar, the nontrivial solutions of these linear equations are the following:
\begin{gather*}
A_{jk}=
\begin{pmatrix}
1 & \pm i
\\
\pm i & -1
\end{pmatrix}
\qquad
\text{with corresponding}
\quad
\lambda=\mp(x_j+x_k),
\\
A_{jk}=
\begin{pmatrix}
1 & \mp i
\\
\pm i & 1
\end{pmatrix}
\qquad
\text{with corresponding}
\quad
\lambda=\mp(x_j-x_k).
\end{gather*}

Let $\epsilon_j\in\mathfrak h^*$ be def\/ined by $\epsilon_j(H)=x_j$ for $1\le j\le\ell$.
Then for $1\le j<k\le\ell$, the following matrices are root vectors of $\mathfrak{so}(2\ell,\mathbb{C})$:
\begin{gather}
X_{\epsilon_j+\epsilon_k} =I_{2k-1,2j-1}-I_{2k,2j}-i(I_{2k-1,2j}+I_{2k,2j-1}),
\nonumber
\\
X_{-\epsilon_j-\epsilon_k} =I_{2k-1,2j-1}-I_{2k,2j}+i(I_{2k-1,2j}+I_{2k,2j-1}),
\nonumber
\\
X_{\epsilon_j-\epsilon_k} =I_{2k-1,2j-1}+I_{2k,2j}-i(I_{2k-1,2j}-I_{2k,2j-1}),
\nonumber
\\
X_{-\epsilon_j+\epsilon_k} =I_{2k-1,2j-1}+I_{2k,2j}+i(I_{2k-1,2j}-I_{2k,2j-1}).
\label{rootvectors}
\end{gather}
Thus, if we choose the following set of positive roots
\begin{gather*}
\Delta^+=\{\epsilon_j+\epsilon_k, \epsilon_j-\epsilon_k: 1\le j<k\le\ell\},
\end{gather*}
then the Dynkin diagram of $\mathfrak{so}(2\ell,\mathbb{C})$ is $D_\ell$:

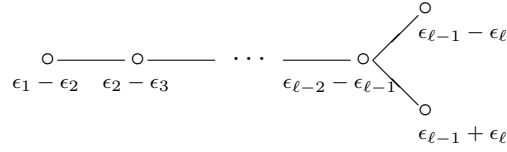
\setlength{\unitlength}{0.75mm} \hspace{-2cm}
\begin{picture}
(0,15) \put(59,0){$\circ$} \put(62,1){\line(1,0){12}} \put(54,-4){\scriptsize$\epsilon_1-\epsilon_2$}
\put(75,0){$\circ$} \put(78,1){\line(1,0){12}} \put(70,-4){\scriptsize$\epsilon_2-\epsilon_3$} \put(93,1){$\dots$}
\put(102,1){\line(1,0){12}}\put(102,-4){\scriptsize$\epsilon_{\ell-2}-\epsilon_{\ell-1}$} \put(115,0){$\circ$}
\put(118,1){\line(1,1){8}} \put(126,9){$\circ$} \put(126,5){\scriptsize$\epsilon_{\ell-1}-\epsilon_{\ell}$}
\put(118,1){\line(1,-1){8}} \put(126,-9){$\circ$} \put(126,-13){\scriptsize$\epsilon_{\ell-1}+\epsilon_{\ell}$}
\end{picture}

\vspace{1cm}

By looking at the $2\times2$ blocks $A_{jk}$ of the dif\/ferent roots, namely
\begin{gather*}
X_{\epsilon_j+\epsilon_k} =
\begin{pmatrix}
1&-i
\\
-i&-1
\end{pmatrix},
\qquad
X_{-\epsilon_j-\epsilon_k}=
\begin{pmatrix}
1&i
\\
i&-1
\end{pmatrix},
\\
X_{\epsilon_j-\epsilon_k}=
\begin{pmatrix}
1&i
\\
-i&1
\end{pmatrix},
\qquad
X_{-\epsilon_j+\epsilon_k}=
\begin{pmatrix}
1&-i
\\
i&1
\end{pmatrix},
\end{gather*}
it is easy to obtain the following inverse relations
\begin{gather*}
I_{2k-1,2j-1} =\tfrac14\big(X_{\epsilon_j+\epsilon_k}+X_{-\epsilon_j-\epsilon_k}+X_{\epsilon_j-\epsilon_k}+X_{-\epsilon_j+\epsilon_k}\big),
\\
I_{2k,2j} =\tfrac14\big(-X_{\epsilon_j+\epsilon_k}-X_{-\epsilon_j-\epsilon_k}+X_{\epsilon_j-\epsilon_k}+X_{-\epsilon_j+\epsilon_k}\big),
\\
I_{2k,2j-1} =\tfrac{i}4\big(X_{\epsilon_j+\epsilon_k}-X_{-\epsilon_j-\epsilon_k}-X_{\epsilon_j-\epsilon_k}+X_{-\epsilon_j+\epsilon_k}\big),
\\
I_{2k-1,2j} =\tfrac{i}4\big(X_{\epsilon_j+\epsilon_k}-X_{-\epsilon_j-\epsilon_k}+X_{\epsilon_j-\epsilon_k}-X_{-\epsilon_j+\epsilon_k}\big).
\end{gather*}
From this it follows that
\begin{gather*}
I_{2k-1,2j-1}^2+I_{2k,2j}^2+I_{2k,2j-1}^2+I_{2k-1,2j}^2
\\
\qquad
=\tfrac14\big(X_{\epsilon_j+\epsilon_k}X_{-\epsilon_j-\epsilon_k}+X_{-\epsilon_j-\epsilon_k}X_{\epsilon_j+\epsilon_k}
+X_{\epsilon_j-\epsilon_k}X_{-\epsilon_j+\epsilon_k}+X_{-\epsilon_j+\epsilon_k}X_{\epsilon_j-\epsilon_k}\big).
\end{gather*}
Therefore
\begin{gather*}
Q_{2\ell}=
\sum\limits_{1\le j\le\ell}I_{2j,2j-1}^2+\tfrac14\sum\limits_{1\le
j<k\le\ell}\big(X_{\epsilon_j+\epsilon_k}X_{-\epsilon_j-\epsilon_k}+X_{-\epsilon_j-\epsilon_k}X_{\epsilon_j+\epsilon_k}
\\
\phantom{Q_{2\ell}=}
{}+X_{\epsilon_j-\epsilon_k}X_{-\epsilon_j+\epsilon_k}+X_{-\epsilon_j+\epsilon_k}X_{\epsilon_j-\epsilon_k}\big).
\end{gather*}
Now using the expressions in~\eqref{rootvectors} we get
\begin{gather*}
[X_{\epsilon_j+\epsilon_k},X_{-\epsilon_j-\epsilon_k}] =-4i(I_{2j,2j-1}+I_{2k,2k-1}),
\\
[X_{\epsilon_j-\epsilon_k},X_{-\epsilon_j+\epsilon_k}] =-4i(I_{2j,2j-1}-I_{2k,2k-1}).
\end{gather*}
Thus $Q_{2\ell}$ becomes
\begin{gather}
Q_{2\ell}=\sum\limits_{1\le j\le\ell}I_{2j,2j-1}^2-2\sum\limits_{1\le j\le\ell}(\ell-j)iI_{2j,2j-1}
\nonumber
\\
\phantom{Q_{2\ell}=}
{}+\sum\limits_{1\le j<k\le\ell}
\tfrac12\big(X_{-\epsilon_j-\epsilon_k}X_{\epsilon_j+\epsilon_k}+X_{-\epsilon_j+\epsilon_k}X_{\epsilon_j-\epsilon_k}\big).
\label{qpar}
\end{gather}

\subsection[The operator $Q_{2\ell+1}$]{The operator $\boldsymbol{Q_{2\ell+1}}$}

Now we look at a~root space decomposition of $\mathfrak{so}(n)$ in terms of the basis elements $I_{ki}$, $1\le i<k\le n$
when $n=2\ell+1$.

The linear span
\begin{gather*}
\mathfrak h={\langle I_{21},I_{43},\dots,I_{2\ell,2\ell-1}\rangle}_\mathbb{C}
\end{gather*}
is a~Cartan subalgebra of $\mathfrak{so}(n,\mathbb{C})$.
To f\/ind the root vectors it is convenient to visualize the elements of $\mathfrak{so}(n,\mathbb{C})$ as $\ell\times\ell$
matrices of $2\times2$ blocks occupying the left upper corner of the square matrices of size $2\ell+1$, with the last
column (respectively row) made up of~$\ell$ columns (respectively rows) of size two and a~zero in the place
$(2\ell+1,2\ell+1)$.
The subspaces of all matrices~$A$ with a~block $A_{jk}$, $1\le j<k\le \ell$, in the place $(j,k)$, with the block
$-A_{jk}^t$ in the place $(k,j)$ and with zeros in all other places, are ${\operatorname{ad}}(\mathfrak h)$-stable.
Also the subspaces of all matrices~$B$ with a~column $B_j$ of size two, $1\le j\le \ell$, in the place $(j,\ell+1)$,
with the row $-B_{j}^t$ in the place $(\ell+1,j)$ and with zeros in all other places, are ${\operatorname{ad}}(\mathfrak
h)$-stable.

On the other hand $[H,B]=\lambda B$ if and only if
\begin{gather*}
x_jiI_{2j,2j-1}B_{j}=\lambda B_{j}.
\end{gather*}
Up to a~scalar this linear equation has two linearly independent solutions:
\begin{gather*}
B_{j}=
\begin{pmatrix}
1
\\
\pm i
\end{pmatrix}
\qquad
\text{with corresponding}
\quad
\lambda=\mp x_j,
\end{gather*}

Let $\epsilon\in\mathfrak h^*$ be def\/ined by $\epsilon(H)=x_j$ for $1\le j\le\ell$.
Then for $1\le j<k\le\ell$ and $1\le r\le\ell$, the following matrices are root vectors of
$\mathfrak{so}(2\ell+1,\mathbb{C})$:
\begin{gather*}
X_{\epsilon_j+\epsilon_k} =I_{2k-1,2j-1}-I_{2k,2j}-i(I_{2k-1,2j}+I_{2k,2j-1}),
\\
X_{-\epsilon_j-\epsilon_k} =I_{2k-1,2j-1}-I_{2k,2j}+i(I_{2k-1,2j}+I_{2k,2j-1}),
\\
X_{\epsilon_j-\epsilon_k} =I_{2k-1,2j-1}+I_{2k,2j}-i(I_{2k-1,2j}-I_{2k,2j-1}),
\\
X_{-\epsilon_j+\epsilon_k} =I_{2k-1,2j-1}+I_{2k,2j}+i(I_{2k-1,2j}-I_{2k,2j-1}),
\\
X_{\epsilon_r} =I_{n,2r-1}-iI_{n,2r},
\\
X_{-\epsilon_r} =I_{n,2r-1}+iI_{n,2r}.
\end{gather*}
Thus, if we choose the following set of positive roots
\begin{gather*}
\Delta^+=\{\epsilon_r, \epsilon_j+\epsilon_k, \epsilon_j-\epsilon_k: 1\le r\le\ell, 1\le j<k\le\ell\},
\end{gather*}
then the Dynkin diagram of $\mathfrak{so}(2\ell+1,\mathbb{C})$ is $B_\ell$:

\setlength{\unitlength}{0.75mm} \hspace{-2cm}
\begin{picture}
(0,15) \put(59,0){$\circ$} \put(62,1){\line(1,0){12}} \put(54,-4){\scriptsize$\epsilon_1-\epsilon_2$}
\put(75,0){$\circ$} \put(78,1){\line(1,0){12}} \put(70,-4){\scriptsize$\epsilon_2-\epsilon_3$} \put(93,1){$\dots$}
\put(102,1){\line(1,0){12}}\put(108,-4){\scriptsize$\epsilon_{\ell-1}-\epsilon_{\ell}$} \put(115,0){$\circ$}
\put(118,1){\line(1,0){12}} \put(118,2){\line(1,0){12}} \put(128,0){$>$} \put(132,0){$\circ$}
\put(131,-4){\scriptsize$\epsilon_{\ell}$}
\end{picture}

\vspace{1cm}

By looking at the $2\times 1$ columns of the dif\/ferent roots, namely
\begin{gather*}
X_{\epsilon_j}=
\begin{pmatrix}
1
\\
-i
\end{pmatrix},
\qquad
X_{-\epsilon_j}=
\begin{pmatrix}
1
\\
i
\end{pmatrix},
\end{gather*}
it is easy to obtain the following inverse relations
\begin{gather*}
I_{n,2r-1}=\tfrac12(X_{\epsilon_r}+X_{-\epsilon_r}),
\qquad
I_{n,2r}=\tfrac{i}2(X_{\epsilon_r}-X_{-\epsilon_r}).
\end{gather*}
From this it follows that
\begin{gather*}
I_{n,2r-1}^2+I_{n,2r}^2 =\tfrac12(X_{\epsilon_r}X_{-\epsilon_r}+X_{-\epsilon_r}X_{\epsilon_r})=-iI_{2r,2r-1}+X_{-\epsilon_r}X_{\epsilon_r},
\end{gather*}
since $[X_{\epsilon_r},X_{-\epsilon_r}]=-2iI_{2r,2r-1}$.
Therefore we have that
\begin{gather*}
Q_{2\ell+1}=\sum\limits_{1\le j\le2\ell}I_{n,j}^2+Q_{2\ell}=\sum\limits_{1\le
r\le2\ell}(-iI_{2r,2r-1}+X_{-\epsilon_r}X_{\epsilon_r})+Q_{2\ell}.
\end{gather*}
Then
\begin{gather}
Q_{2\ell+1}=\sum\limits_{1\le j\le\ell}I_{2j,2j-1}^2-\sum\limits_{1\le j\le\ell}(2(\ell-j)+1)iI_{2j,2j-1}
\nonumber
\\
\phantom{Q_{2\ell+1}=}
{}+\sum\limits_{1\le j<k\le\ell}
\tfrac12\big(X_{-\epsilon_j-\epsilon_k}X_{\epsilon_j+\epsilon_k}+X_{-\epsilon_j+\epsilon_k}X_{\epsilon_j-\epsilon_k}\big)
+\sum\limits_{1\le r\le2\ell}X_{-\epsilon_r}X_{\epsilon_r}.
\label{qimpar}
\end{gather}

\subsection{Gel'fand--Tsetlin basis}
\label{GTB}

For any~$n$ we identify the group ${\mathrm{SO}}(n)$ with a~subgroup of ${\mathrm{SO}}(n+1)$ in the following way: given
$k\in{\mathrm{SO}}(n)$ we have
\begin{gather*}
k\simeq
\begin{pmatrix}
k&\bf0
\\
\bf0&1
\end{pmatrix}
\in{\mathrm{SO}}(n+1).
\end{gather*}

Let $T_{\mathbf{m}}$ be an irreducible unitary representation of ${\mathrm{SO}}(n)$ with highest weight ${\mathbf{m}}$
and let $V_{\mathbf{m}}$ be the space of this representation.
Highest weights ${\mathbf{m}}$ of these representations are given by the~$\ell$-tuples of integers
${\mathbf{m}}={\mathbf{m}}_n=(m_{1n},\dots,m_{\ell n})$ for which
\begin{gather*}
m_{1n}\ge m_{2n}\ge \dots\ge m_{\ell-1,n}\ge|m_{\ell n}|
\qquad
\text{if}
\quad
n=2\ell,
\\
m_{1n}\ge m_{2n}\ge \dots\ge m_{\ell n}\ge0
\qquad
\text{if}
\quad
n=2\ell+1,
\end{gather*}
and $m_{jn}$ are all integers.

The restriction of the representation $T_{\mathbf{m}}$ of the group ${\mathrm{SO}}(2\ell+1)$ to the subgroup
${\mathrm{SO}}(2\ell)$ decomposes into the direct sum of all representations $T_{{\mathbf{m}}'}$,
${\mathbf{m}}'={\mathbf{m}}_{n-1}=(m_{1,n-1},\dots,m_{\ell,n-1})$ for which the betweenness conditions
\begin{gather*}
m_{1,2\ell+1}\ge m_{1,2\ell}\ge m_{2,2\ell+1}\ge m_{2,2\ell}\ge\dots\ge m_{\ell,2\ell+1}\ge
m_{\ell,2\ell}\ge-m_{\ell,2\ell+1}
\end{gather*}
are satisf\/ied.
For the restriction of the representations $T_{\mathbf{m}}$ of ${\mathrm{SO}}(2\ell)$ to the subgroup
${\mathrm{SO}}(2\ell-1)$ the corresponding betweenness conditions are
\begin{gather*}
m_{1,2\ell}\ge m_{1,2\ell-1}\ge m_{2,2\ell}\ge m_{2,2\ell-1}\ge\dots\ge m_{\ell-1,2\ell}\ge
m_{\ell-1,2\ell-1}\ge|m_{\ell,2\ell}|.
\end{gather*}
All multiplicities in the decompositions are equal to one (see~\cite[p.~362]{V92}).

If we continue this procedure of restriction of irreducible representations successively to the subgroups
\begin{gather*}
{\mathrm{SO}}(n-2)>{\mathrm{SO}}(n-3)>\dots>{\mathrm{SO}}(2),
\end{gather*}
then we f\/inally get one dimensional representations of the group ${\mathrm{SO}}(2)$.
If we take a~unit vector in each one of these one dimensional representations we get an orthonormal basis of the
representation space $V_{\mathbf{m}}$.
Such a~basis is called a~Gel'fand--Tsetlin basis.
The elements of a~Gel'fand--Tsetlin basis $\{v(\mu)\}$ of the representation $T_{\mathbf{m}}$ of ${\mathrm{SO}}(n)$ are
labelled by the Gel'fand--Tsetlin patterns
$\mu=({\mathbf{m}}_{n},{\mathbf{m}}_{n-1},\dots,{\mathbf{m}}_3,{\mathbf{m}}_2)$, where the betweenness conditions are
depicted in the following diagrams.

If $n=2\ell+1$
\[
\mu=
\begin{tikzpicture}[baseline=-0.5ex]\small
\matrix (magic) [matrix of math nodes, left delimiter = {[}, right delimiter = {]}]
  {%
 m_{1n}&{}&m_{2n}&{}&{}&{}&{}&m_{\ell\, n}&{}&-m_{\ell\, n} \\
{}&m_{1,n-1}&{}&{}&{}&{}&{}&{}&m_{\ell,n-1}&{} \\
{\phantom{m_{2n}}}&{}&{}&{}&{}&{}&{}&{}&{}&{}\\
{}&{}&{}&{}&{}&m_{15}&{}&m_{25}&{}&-m_{25}\\
{}&{}&{}&{}&{}&{}&m_{14}&{}&m_{24}&{}\\
{}&{}&{}&{}&{}&{}&{}&m_{13}&{}&-m_{13}\\
{}&{}&{}&{}&{}&{}&{}&{}&m_{12}&{}\\
  };
  \draw[loosely dotted,line width=1] (magic-1-5) -- (magic-1-7
  );
  \draw[loosely dotted,line width=1] (magic-2-3) -- (magic-2-8
  );
  \draw[loosely dotted,line width=1] (magic-3-3.west) -- (magic-3-10.east);
\end{tikzpicture}.
\]

If $n=2\ell$
\[
\mu=
\begin{tikzpicture}[baseline=-0.5ex]\small
\matrix (magic) [matrix of math nodes, left delimiter = {[}, right delimiter = {]}]
  {%
 m_{1n}&{}&m_{2n}&{}&{}&{}&{}&m_{\ell\, n}&{}\\
{}&m_{1,n-1}&{}&{}&{}&{}&m_{\ell-1,n-1}&{}&-m_{\ell-1,n-1}\\
{\phantom{m_{2n}}}&{}&{}&{}&{}&{}&{}&{}&{}\\
{}&{}&{}&{}&m_{15}&{}&m_{25}&{}&-m_{25}\\
{}&{}&{}&{}&{}&m_{14}&{}&m_{24}&{}\\
{}&{}&{}&{}&{}&{}&m_{13}&{}&-m_{13}\\
{}&{}&{}&{}&{}&{}&{}&m_{12}&{}\\
  };
  \draw[loosely dotted,line width=1] (magic-1-5) -- (magic-1-7
  );
  \draw[loosely dotted,line width=1] (magic-2-3) -- (magic-2-6
  );
  \draw[loosely dotted,line width=1] (magic-3-3.west) -- (magic-3-10.east);
\end{tikzpicture}.
\]

The chain of subgroups ${\mathrm{SO}}(n-1)>{\mathrm{SO}}(n-2)>\dots>{\mathrm{SO}}(2)$ def\/ines the orthonormal basis
$\{v(\mu)\}$ uniquely up to multiplication of the basis elements by complex numbers of absolute value one.

\subsection[An explicit expression for $\dot\pi(Q_n)$]{An explicit expression for $\boldsymbol{\dot\pi(Q_n)}$}

Since $Q_n\in D({\mathrm{SO}}(n))^{{\mathrm{SO}}(n)}$, given $\dot\pi\in\hat{\mathrm{SO}}(n)$ it follows that
$\dot\pi(Q_n)$ commutes with $\pi(k)$ for all $k\in{\mathrm{SO}}(n)$.
Hence, by Schur's Lemma $\dot\pi(Q_n)=\lambda I$.
From expressions~\eqref{qpar} and~\eqref{qimpar} we can give the explicit value of~$\lambda$ in terms of the highest
weight of~$\pi$, by computing $\dot\pi(Q_n)$ on a~highest weight vector.

\begin{Proposition}
Let $(\pi,V_\pi)$ be an irreducible representation of ${\mathrm{SO}}(2\ell)$ of highest weight
${\mathbf{m}}=(m_1,m_2,\dots,m_\ell)$.
Then, $\dot\pi(Q_{2\ell})=\lambda I$, with
\begin{gather}
\label{qvpar}
\lambda= \sum\limits_{1\le j\le\ell}\big({-}m_j^2-2(\ell-j)m_j\big).
\end{gather}
\end{Proposition}

\begin{Proposition}
Let $(\pi,V_\pi)$ be an irreducible representation of ${\mathrm{SO}}(2\ell+1)$ of highest weight
${\mathbf{m}}=(m_1,m_2,\dots,m_\ell)$.
Then, $\dot\pi(Q_{2\ell+1})=\lambda I$, with
\begin{gather}
\label{qvimpar}
\lambda= \sum\limits_{1\le j\le\ell}\left(-m_j^2-(2(\ell-j)+1)m_j\right).
\end{gather}
\end{Proposition}

\section[The dif\/ferential operator~$\Delta$]{The dif\/ferential operator~$\boldsymbol{\Delta}$}
\label{sec: delta}

We shall look closely at the left invariant dif\/ferential operator~$\Delta$ of ${\mathrm{SO}}(n+1)$ def\/ined~by
\begin{gather*}
\Delta=\sum\limits_{j=1}^{n} I_{n+1,j}^2,
\end{gather*}
in order to study its eigenfunctions and eigenvalues.
Later we will use all this to understand the irreducible spherical functions of fundamental $K$-types associated with the
pair $(G,K)=({\mathrm{SO}}(n+1),{\mathrm{SO}}(n))$.

\begin{Proposition}
Let $G={\mathrm{SO}}(n+1)$ and $K={\mathrm{SO}}(n)$.
Let us consider the following left invariant differential operator of~$G$
\begin{gather*}
\Delta=\sum\limits_{j=1}^{n} I_{n+1,j}^2.
\end{gather*}
Then~$\Delta$ is also right invariant under~$K$.
\end{Proposition}
\begin{proof}
This is a~direct consequence of the identity
\begin{gather*}
Q_{n+1}=Q_{n}+\Delta
\end{gather*}
and Proposition~\ref{rightinv}.
\end{proof}

Let us def\/ine the one-parameter subgroup~$A$ of~$G$ as the set of all elements of the form
\begin{gather}
\label{a}
a(s)=
\begin{pmatrix}
I_{n-1}&{\bf 0}&{\bf 0}
\\
{\bf 0}&\cos s&\sin s
\\
{\bf 0} &-\sin s&\cos s
\end{pmatrix},
\qquad
-\pi\le s\le\pi,
\end{gather}
where $I_{n-1}$ denotes the identity matrix of size $n-1$, and let $M={\mathrm{SO}}(n-1)$ be the centralizer of~$A$
in~$K$.

Now we want to get the expression of $[\Delta\Phi](a(s))$ for any smooth function~$\Phi$ on~$G$ with values in
${\operatorname{End}}(V_\pi)$ such that $\Phi(kgk')=\pi(k)\Phi(g)\pi(k')$ for all $g\in G$ and all $k,k'\in K$.

We have
\begin{gather*}
\big[I_{n+1,j}^2\Phi\big](a(s))
=\left.\frac{\partial^2}{\partial t^2}\Phi(a(s)\exp tI_{n+1,j})\right|_{t=0}.
\end{gather*}
Hence, we will use the decomposition $G=KAK$ to write $a(s)\exp tI_{n+1,j}=k(s,t)a(s,t)h(s,t)$, with $k(s,t),h(s,t)\in
K$ and $a(s,t)\in A$.

Let us take on $A \setminus \{a(\pi)\}$ the coordinate function $x(a(s))=s$, with $-\pi< s<\pi$, and let
\begin{gather*}
F(s)=F(x(a(s)))=\Phi(a(s)).
\end{gather*}
From now on we will assume that $-\pi<s,t,s+t<\pi$.

If $j=n$ we have $a(s)\exp tI_{n+1,n}=a(s)a(t)=a(s+t)$.
Thus we may take
\begin{gather*}
a(s,t)=a(s+t),
\qquad
k(s,t)=h(s,t)=e.
\end{gather*}

Since $x(a(s+t))=s+t$, we obtain
\begin{gather*}
\big[I_{n+1,n}^2\Phi\big](a(s))=\left.\frac{\partial^2}{\partial t^2}\Phi(a(s)\exp tI_{n+1,n})\right|_{t=0}
=\left.\frac{\partial^2}{\partial t^2}\Phi(a(s+t))\right|_{t=0}
\\
\phantom{\big[I_{n+1,n}^2\Phi\big](a(s))}
 =\left.\frac{\partial^2}{\partial t^2}F(s+t)\right|_{t=0}=F''(s).
\end{gather*}

For $1\le j\le n-1$, when $s\notin\mathbb{Z}\pi$, we may take
\begin{gather*}
k(s,t)=
\begin{pmatrix}
I_{j-1} & {\bf 0} & {\bf 0} & {\bf 0} & {\bf 0}
\\
{\bf 0} & \frac{\sin s\cos t}{\sqrt{1-\cos^2s\cos^2t}} & {\bf 0} & \frac{\sin t}{\sqrt{1-\cos^2s\cos^2t}} & 0
\\
{\bf 0} & {\bf 0} & I_{n-j-1} & {\bf 0} & {\bf 0}
\\
{\bf 0} & \frac{-\sin t}{\sqrt{1-\cos^2s\cos^2t}} & {\bf 0} & \frac{\sin s\cos t}{\sqrt{1-\cos^2s\cos^2t}} & 0
\\
{\bf 0} & 0 & {\bf 0} & 0 & 1
\end{pmatrix},
\\
h(s,t)=
\begin{pmatrix}
I_{j-1} & {\bf 0} & {\bf 0} & {\bf 0} & {\bf 0}
\\
{\bf 0} & \frac{\sin s}{\sqrt{1-\cos^2s\cos^2t}} & {\bf 0} & \frac{-\cos s\sin t}{\sqrt{1-\cos^2s\cos^2t}} & 0
\\
{\bf 0} & {\bf 0} & I_{n-j-1} & {\bf 0} & {\bf 0}
\\
{\bf 0} & \frac{\cos s\sin t}{\sqrt{1-\cos^2s\cos^2t}} & {\bf 0} & \frac{\sin s}{\sqrt{1-\cos^2s\cos^2t}} & 0
\\
{\bf 0} & 0 & {\bf 0} & 0 & 1
\end{pmatrix},
\\
a(s,t)=
\begin{pmatrix}
I_{n-1} & {\bf 0} & {\bf 0}
\\
{\bf 0} & \cos s\cos t & \sqrt{1-\cos^2s\cos^2t}
\\
{\bf 0} & -\sqrt{1-\cos^2s\cos^2t} & \cos s\cos t
\end{pmatrix}.
\end{gather*}

Then, for $0<s<\pi$, we have $x(a(s,t))=\arccos(\cos s\cos t)$ and
\begin{gather*}
\frac{\partial}{\partial t}x(a(s,t))=\frac{\cos s\sin t}{\sqrt{1-\cos^2 s\cos^2t}}.
\end{gather*}
From here we get
\begin{gather*}
\left.\frac{\partial}{\partial t}x(a(s,t))\right|_{t=0}=0
\qquad
\text{and}
\qquad
\left.\frac{\partial^2}{\partial t^2}x(a(s,t))\right|_{t=0}=\frac{\cos s}{\sin s}.
\end{gather*}
Thus
\begin{gather*}
\left.\frac{\partial}{\partial t}\Phi(a(s,t))\right|_{t=0}=\left.F'(s)\frac{\partial}{\partial t}x(a(s,t))\right|_{t=0}=0
\qquad
\text{and}
\qquad
\left.\frac{\partial^2}{\partial t^2}\Phi(a(s,t))\right|_{t=0}=\frac{\cos s}{\sin s}F'(s).
\end{gather*}

We observe that $k(s,0)=h(s,0)=e$ and that $a(s,0)=a(s)$.
Then
\begin{gather*}
[I_{nj}^2\Phi](a(s))= \frac{\partial^2}{\partial t^2}\pi(k(s,t))\Big|_{t=0}\Phi(a(s))+2\frac{\partial}{\partial
t}\pi(k(s,t))\Big|_{t=0}\frac{\partial}{\partial t}\Phi(a(s,t))\Big|_{t=0}
\\
\phantom{[I_{nj}^2\Phi](a(s))=}
{}+2\frac{\partial}{\partial t}\pi(k(s,t))\Big|_{t=o}\Phi(a(s))\frac{\partial}{\partial
t}\pi(h(s,t))\Big|_{t=0}+\frac{\partial^2}{\partial t^2}\Phi(a(s,t))\Big|_{t=0}
\\
\phantom{[I_{nj}^2\Phi](a(s))=}
{}+2\frac{\partial}{\partial t}\Phi(a(s,t))\Big|_{t=0}\frac{\partial}{\partial
t}\pi(h(s,t))\Big|_{t=0}+\Phi(a(s))\frac{\partial^2}{\partial t^2}\pi(h(s,t))\Big|_{t=0}.
\end{gather*}

We also have
\begin{gather*}
\frac{\partial}{\partial t}\pi\left(k(s,t)\right)\Big|_{t=0}=\dot\pi\left(\frac{\partial}{\partial
t}k(s,t)\Big|_{t=0}\right)=\frac1{\sin s}\dot\pi(I_{n,j}),
\end{gather*}
and
\begin{gather*}
\frac{\partial}{\partial t}\pi(h(s,t))\Big|_{t=0}=\dot\pi\left(\frac{\partial}{\partial
t}h(s,t)\Big|_{t=0}\right)=-\frac{\cos s}{\sin s}\dot\pi\left(I_{n,j}\right).
\end{gather*}

We will need the following proposition, whose proof appears in the Appendix and its idea is taken from \cite{GPT02a}.
\begin{Proposition}
\label{derprincipal}
If $A(s,t)=k(s,t)$ or $A(s,t)=h(s,t)$, then in either case for $0<s<\pi$, we have
\begin{gather*}
\frac{\partial^2(\pi\circ A)}{\partial t^2}\Big|_{t=0}=\dot\pi\left(\frac{\partial A}{\partial t}\Big|_{t=0}\right)^2.
\end{gather*}
Moreover in each case, for $1\le j\le n-1$ and $0<s<\pi$, we have
\begin{gather*}
\frac{\partial^2}{\partial t^2}\pi(k(s,t))\Big|_{t=0}=\frac1{\sin^2s}\dot\pi(I_{n,j})^2,
\qquad
\frac{\partial^2}{\partial t^2}\pi(h(s,t))\Big|_{t=0}=\frac{\cos^2 s}{\sin^2 s}\dot\pi(I_{n,j})^2.
\end{gather*}
\end{Proposition}

Now we obtain the following corollaries.
\begin{Corollary}
Let~$\Phi$ be any smooth function on~$G$ with values in ${\operatorname{End}}(V_\pi)$ such that
$\Phi(kgk')=\pi(k)\Phi(g)\pi(k')$ for all $g\in G$ and all $k,k'\in K$.
Then, if $F(s)=\Phi(a(s))$, for $0<s<\pi$ we have
\begin{gather*}
[\Delta\Phi](a(s))=F''(s)+(n-1)\frac{\cos s}{\sin s}F'(s) +\frac{1}{\sin^2
s}\sum\limits_{j=1}^{n-1}\dot\pi(I_{n,j})^2F(s)
\\
\phantom{[\Delta\Phi](a(s))=}
{}-2\frac{\cos s}{\sin^2 s}\sum\limits_{j=1}^{n-1}\dot\pi(I_{n,j})F(s)\dot\pi(I_{n,j}) +\frac{\cos^2s}{\sin^2
s}F(s)\sum\limits_{j=1}^{n-1}\dot\pi(I_{n,j})^2.
\end{gather*}
\end{Corollary}

\begin{Corollary}
\label{eigenfunction}
Let~$\Phi$ be an irreducible spherical function on~$G$ of type $\pi\in \hat K$.
Then, if $F(s)=\Phi(a(s))$, we have
\begin{gather*}
F''(s)+(n-1)\frac{\cos s}{\sin s}F'(s)+\frac1{\sin^2s}\sum\limits_{j=1}^{n-1}\dot\pi(I_{n,j})^2F(s)
\\
\qquad
{}-2\frac{\cos s}{\sin^2
s}\sum\limits_{j=1}^{n-1}\dot\pi(I_{n,j})F(s)\dot\pi(I_{n,j})+\frac{\cos^2s}{\sin^2s}F(s)\sum\limits_{j=1}^{n-1}\dot\pi(I_{n,j})^2=\lambda
F(s),
\end{gather*}
for some $\lambda\in\mathbb{C}$ and $0<s<\pi$.
\end{Corollary}

Notice that the expression in Corollary~\ref{eigenfunction} generalizes the very well known situation when the $K$-type
is the trivial one, as we state in the following corollary (cf.~\cite[p.~403, equation~(10)]{Helgason1962}).
\begin{Corollary}
Let~$\Phi$ be an irreducible spherical function on~$G$ of the trivial $K$-type.
Then, for $F(s)=\Phi(a(s))$ we have
\begin{gather*}
 F''(s)+(n-1)\frac{\cos s}{\sin s}F'(s)=\lambda F(s),
\end{gather*}
for some $\lambda\in\mathbb{C}$ and $0<s<\pi$.
\end{Corollary}

Let us make the change of variables $y=(1+\cos s)/2$, with $0<s<\pi$; then $0<y<1$.
We also have $\cos s=2y-1$, $\sin^2s=4y(1-y)$ and $\tfrac{d}{dy}=-\tfrac{\sin s}{2}$.
If we let $H(y)=F(s)$, i.e.~
\begin{gather*}
H(y)=\Phi(a(s)),
\qquad
\text{with}
\quad
\cos s=2y-1,
\end{gather*}
we obtain
\begin{gather*}
F'(s)=-\frac{\sin s}{2}H'(s),
\qquad
F''(s)=\frac{\sin^2s}{4}H''(y)-\frac{\cos s}{2}H'(y).
\end{gather*}

In terms of this new variable Corollary~\ref{eigenfunction} becomes
\begin{Corollary}
\label{eigenfunction2}
Let~$\Phi$ be an irreducible spherical function on~$G$ of type $\pi\in \hat K$.
Then, if $H(y)=\Phi(a(s))$ with $y=(1+\cos s)/2$, we have
\begin{gather*}
y(1-y)H''(y)+\frac{1}{2}n(1-2y)H'(y)+\frac1{4y(1-y)}\sum\limits_{j=1}^{n-1}\dot\pi(I_{n,j})^2H(y)
\\
\qquad
{}+\frac{(1-2y)}{2y(1-y)}\sum\limits_{j=1}^{n-1}\dot\pi(I_{n,j})H(y)\dot\pi(I_{n,j})
+\frac{(1-2y)^2}{4y(1-y)}H(y)\sum\limits_{j=1}^{n-1}\dot\pi(I_{n,j})^2=\lambda H(y),
\end{gather*}
for some $\lambda\in\mathbb{C}$ and $0<y<1$.
\end{Corollary}

\begin{Remark}
\label{hache}
Let us notice that, for any $y\in (0,1)$, $H(y)$ is a~scalar linear transformation when restricted to any~$M$-submodule,
see Proposition~\ref{propesf}.
Therefore, if~$m$ is the number of~$M$-submodules contained in $(V,\pi)$, we consider the vector valued function
$H:(0,1)\to \mathbb{C}^m$ whose entries are given by those scalar values that $H(y)$ takes on every~$M$-submodule.
\end{Remark}

If the $\operatorname{End}(V)$-valued function~$H$ satisf\/ies the dif\/ferential equation given in
Corollary~\ref{eigenfunction2}, then the vector valued function~$H$ satisf\/ies
\begin{gather*}
y(1-y)H''(y)+\frac{1}{2}n(1-2y)H'(y)+\frac1{4y(1-y)}N_1H(y)
\\
\qquad
{}+\frac{(1-2y)}{2y(1-y)}E H(y)+\frac{(1-2y)^2}{4y(1-y)}N_2 H(y)=\lambda H(y),
\end{gather*}
where~$E$, $N_1$ and $N_2$ are matrices of size $m\times m$.

Even more, since $\sum\limits_{j=1}^{n-1}I_{n,j}^2=Q_n-Q_{n-1}$, Proposition~\ref{rightinv} implies
$\sum\limits_{j=1}^{n-1}I_{n,j}^2\in D({\mathrm{SO}}(n))^{{\mathrm{SO}}(n-1)}$, therefore
$\sum\limits_{j=1}^{n-1}\dot\pi(I_{n,j})^2$ is scalar valued when restricted to any~$M$-sub\-mo\-du\-le.
Hence, $N_1=N_2$ and the equation above is equivalent to
\begin{gather}
\label{vectoreq}
y(1-y)H''(y) \! + \! \frac{n}{2}(1-2y)H'(y)+\frac{(1-2y)}{2y(1-y)}E H(y)+\frac{1+(1-2y)^2}{4y(1-y)}N H(y) \! = \! \lambda H(y),  \!\!\!
\end{gather}
where~$N$ is a~diagonal matrix of size $m\times m$.
To obtain an explicit expression of~$E$ for any $K$-type is a~very serious matter; in the following sections we shall
f\/ind explicitly the expressions of~$E$ and~$N$, for certain $K$-types.
\begin{Remark}
\label{ene}
It is worth to observe that from~\eqref{qvpar} and~\eqref{qvimpar} we can immediately obtain every entry of the diagonal
matrix~$N$.
\end{Remark}

\section[The $K$-types which are~$M$-irreducible]{The $\boldsymbol{K}$-types which are~$\boldsymbol{M}$-irreducible}
\label{mirred}

Let $K={\mathrm{SO}}(n)$, $M={\mathrm{SO}}(n-1)$, with $n=2\ell+1$, and let ${\mathbf{m}}_n=(m_{1n},\dots,m_{\ell n})$
be a $K$-type such that $V_{\mathbf{m}}$ is irreducible as~$M$-module.
The highest weights ${\mathbf{m}}_{n-1}$ of the~$M$-submodules of~$V_{\mathbf{m}}$ are those that satisf\/ies the
following intertwining relations
\begin{gather*}
\begin{matrix} m_{1n}&{} &m_{2n} &{} &\dots&{} &m_{\ell,n}&{} &-m_{\ell n}
\\
{}&m_{1,n-1}&{} &\dots&{} &\dots &{} &m_{\ell,n-1}&{}
\end{matrix}.
\end{gather*}

Since $V_{\mathbf{m}}$ is irreducible as~$M$-module it follows that $m_{1n}=\dots=m_{\ell,n}=0$.
The converse is also true, therefore $V_{\mathbf{m}}$ is~$M$-irreducible if and only if it is the trivial
representation.

Let now consider the case $K={\mathrm{SO}}(n)$, $M={\mathrm{SO}}(n-1)$, with $n=2\ell$ and let
${\mathbf{m}}_n=(m_{1n},\dots,m_{\ell n})$ be a $K$-type such that $V_{\mathbf{m}}$ is irreducible as~$M$-module.
The highest weights ${\mathbf{m}}_{n-1}$ of the~$M$-submodules of $V_{\mathbf{m}}$ are those that satisf\/ies the
following intertwining relations
\begin{gather*}
\begin{matrix}
m_{1n}&{}&m_{2n}&{}&\dots&{}&m_{\ell-1,n}&{}&m_{\ell n}&{}
\\
{}&m_{1,n-1}&{}&\dots&{}&\dots&{}&m_{\ell-1,n-1}&{}&-m_{\ell-1,n-1}
\end{matrix}
.
\end{gather*}
Since $V_{\mathbf{m}}$ is irreducible as~$M$-module it follows that $m_{1n}=\dots=m_{\ell-1,n}=d$ and $m_{\ell n}=d-j$
with $0\le j\le 2d$, since $m_{\ell-1,n}\ge\vert m_{\ell n}\vert$.
This implies that $m_{1,n-1}=\dots =m_{\ell-2,n-1}=d$ and $m_{\ell-1,n-1}=q$ with $d\ge q\ge\max\{d-j,j-d\}$.
Thus, if $0\le j\le d$ we have $d\ge q\ge d-j$ and by irreducibility we must have $j=0$.
Similarly if $d\le j\le 2d$ we have $d\ge q\ge j-d$ and by irreducibility we must have $j=2d$.
Therefore ${\mathbf{m}}_n=d\alpha$ or ${\mathbf{m}}_n=d\beta$, where
\begin{gather*}
\alpha=(1,\dots,1),
\qquad
\beta=(1,\dots,1,-1).
\end{gather*}

The converse is also true, therefore $V_{\mathbf{m}}$ is~$M$-irreducible if and only if ${\mathbf{m}}_n=d\alpha$ or
${\mathbf{m}}_n=d\beta$ for any $d\in\mathbb{N}_0$.

If~$\Phi$ is an irreducible spherical function on ${\mathrm{SO}}(n+1)$ of type~$\pi$, whose highest weight is
${\mathbf{m}}_n=d\alpha$ or ${\mathbf{m}}_n=d\beta$, then from Corollary~\ref{eigenfunction2} we get that the associated
function~$H$ satisf\/ies
\begin{gather*}
y(1-y)H''(y)+\ell(1-2y)H'(y)+\frac{1-y}{y}\sum\limits_{j=1}^{n-1}\dot\pi(I_{nj})^2H(y)=\lambda H(y).
\end{gather*}

To compute $\sum\limits_{j=1}^{n-1}\dot\pi(I_{nj})^2$ we write
$\sum\limits_{j=1}^{n-1}\dot\pi(I_{nj})^2=\dot\pi(Q_{n}-Q_{n-1})$.

Let us f\/irst consider ${\mathbf{m}}_n=d\alpha$.
If $v\in V_{{\mathbf{m}}_n}$ is a~highest weight vector, then
\begin{gather*}
\dot\pi(Q_n)v=-d\ell(d+\ell-1)v
\qquad
\text{and}
\qquad
\dot\pi(Q_{n-1})v=-d(\ell-1)(d+\ell-1)v,
\end{gather*}
see~\eqref{qvpar} and~\eqref{qvimpar}.
Therefore
\begin{gather*}
\sum\limits_{j=1}^{n-1}\dot\pi(I_{nj})^2v=-d(d+\ell-1)v.
\end{gather*}

Let us now consider ${\mathbf{m}}_n=d\beta$.
If $v\in V_{{\mathbf{m}}_n}$ is a~highest weight vector, then $\dot\pi(Q_n)v=-2d\ell(d+\ell-1)v$ as before, and
$\dot\pi(Q_{n-1})v=-2d(\ell-1)(d+\ell-1)v$ as before because in both cases ${\mathbf{m}}_{n-1}$ is the same.

Therefore if ${\mathbf{m}}_n=(d,\dots,d,\pm d)$ we have
\begin{gather*}
\sum\limits_{j=1}^{n-1}\dot\pi(I_{nj})^2v=-d(d+\ell-1)v.
\end{gather*}

Hence, if~$\Phi$ is an irreducible spherical function on ${\mathrm{SO}}(n+1)$, $n=2\ell$, of type
${\mathbf{m}}_n=(d,\dots,d,\pm d)\in \mathbb{C}^\ell$, then the associated scalar value function $H=h$ satisf\/ies
\begin{gather}
\label{ecuacion}
y(1-y)h''(y)+\ell(1-2y)h'(y)-\frac{d(d+\ell-1)(1-y)}{y}h(y)=\lambda h(y).
\end{gather}

Let us now compute the eigenvalue~$\lambda$ corresponding to the spherical function of type
$\pi\in\hat{\mathrm{SO}}(2\ell)$, of highest weight ${\mathbf{m}}_n=d\alpha$, associated with the irreducible
representation $\tau\in{\mathrm{SO}}(2\ell+1)$, of highest weight ${\mathbf{m}}_{n+1}=(w,d,\dots,d)\in\mathbb{C}^\ell$.
If $v\in V_{{\mathbf{m}}_{n+1}}$ is a~highest weight vector, then from~\eqref{qvimpar} we have
\begin{gather*}
\dot\tau(Q_{n+1})v=-\left(w(w+2\ell-1)+d(\ell-1)(d+\ell-1)\right)v.
\end{gather*}

If $v\in V_{{\mathbf{m}}_{n}}$ is a~highest weight vector, then from~\eqref{qvpar} we have
\begin{gather*}
\dot\tau(Q_n)v=\dot\pi(Q_n)v=-d\ell(d+\ell-1)v.
\end{gather*}
Since $\Delta=Q_{n+1}-Q_n $ it follows that
\begin{gather*}
\lambda=-w(w+2\ell-1)+d(d+\ell-1).
\end{gather*}

To solve~\eqref{ecuacion} we write $h=y^\alpha f$.
Then we get
\begin{gather*}
y(1-y)y^\alpha f''+(2\alpha(1-y)+\ell(1-2y))y^\alpha f'
\\
\qquad
{}+(\alpha(\alpha-1)(1-y)+\ell\alpha(1-2y)-d(d+\ell-1)(1-y))y^{\alpha-1}f=\lambda y^\alpha f.
\end{gather*}
Thus the indicial equation is $\alpha(\alpha-1)+\ell\alpha-d(d+\ell-1)=0$ and $\alpha=d$ is one of its solutions.
If we take $h=y^df$, then we obtain
\begin{gather*}
y(1-y)f''+(2d+\ell-2(d+\ell)y)f'-d\ell f=\lambda f.
\end{gather*}
If we replace $\lambda=-w(w+2\ell-1)+d(d+\ell-1)$ we get
\begin{gather*}
y(1-y)f''+(2d+\ell-2(d+\ell)y)f'-(d-w)(2\ell+d+w-1)f=0.
\end{gather*}

Let $a=d-w$, $b=2\ell+d+w-1$, $c=2d+\ell$ then the above equation becomes
\begin{gather*}
y(1-y)f''+(c-(1+a+b)y)f'-abf=0.
\end{gather*}

A fundamental system of solutions of this equation near $y=0$ is given by the following functions
\begin{gather*}
{}_2F_1 \left(
\begin{matrix}
a,b
\\
c
\end{matrix}
; y\right),
\qquad
y^{1-c}{}_2F_1 \left(
\begin{matrix}
a-c+1,b-c+1
\\
2-c
\end{matrix}
; y\right).
\end{gather*}
Since $h=y^df$ is bounded near $y=0$ it follows that
\begin{gather*}
h(y)=uy^d{}_2F_1 \left(
\begin{matrix}
d-w,2\ell+d+w-1
\\
2d+\ell
\end{matrix};
y\right),
\end{gather*}
where the constant~$u$ is determined by the condition $h(1)=1$.
\begin{Remark}
Let $h_w=h_w(y)$, $w\ge d$, be the function~$h$ above.
Then $h_w$ is a~polynomial of degree~$w$.
Moreover observe that the function $y^d$ used to hypergeometrize~\eqref{ecuacion} is precisely $h_d$.
\end{Remark}

Let us now compute the eigenvalue~$\lambda$ corresponding to the spherical function of type ${\mathbf{m}}_n=d\beta$
associated with an irreducible representation~$\tau$ of ${\mathrm{SO}}(n+1)$ of highest weight
${\mathbf{m}}_{n+1}=(w,d,\dots,d)\in\mathbb{C}^\ell$.
If $v\in V_{{\mathbf{m}}_{n+1}}$ is a~highest weight vector, we obtain $\dot\tau
(Q_{n+1})v=-(w(w+2\ell-1)+d(\ell-1)(d+\ell-1))v$.

If $v\in V_{{\mathbf{m}}_{n}}$ is a~highest weight vector, then $\dot\pi(Q_n)v=-d\ell(d+\ell-1)v$ as above, because
$Q_nv$ does not depend on the sign of the last coordinate of ${\mathbf{m}}_n$.
Since $\Delta=Q_{n+1}-Q_n$ we also have
\begin{gather*}
\lambda=-w(w+2\ell-1)+d(d+\ell-1).
\end{gather*}

Therefore we have proved the following result.
\begin{Theorem}
The scalar valued functions $H=h$ associated with the irreducible spherical functions on ${\mathrm{SO}}(n+1)$,
$n=2\ell$, of ${\mathrm{SO}}(n)$-type ${\mathbf{m}}_n=(d,\dots,d,\pm d)\in \mathbb{C}^\ell$, are parameterized by the
integers $w\ge d$ and are given~by
\begin{gather*}
h_w(y)=uy^d{}_2F_1 \left(
\begin{matrix}
d-w,2\ell+d+w-1
\\
2d+\ell
\end{matrix}
; y\right)
\end{gather*}
where the constant~$u$ is determined by the condition $h_w(1)=1$.
\end{Theorem}

\section[The operator~$\Delta$ for fundamental $K$-types]{The operator~$\boldsymbol{\Delta}$ for fundamental $\boldsymbol{K}$-types}
\label{sec: fund}

We are interested in f\/inding a~more explicit expression of the dif\/ferential equation given in
Corollary~\ref{eigenfunction2}:
\begin{gather*}
y(1-y)H''(y)+\frac{1}{2}n(1-2y)H'(y)+\frac1{4y(1-y)}\sum\limits_{j=1}^{n-1}\dot\pi(I_{n,j})^2H(y)
\\
\qquad
{}+\frac{(1-2y)}{2y(1-y)}\sum\limits_{j=1}^{n-1}\dot\pi(I_{n,j})H(y)\dot\pi(I_{n,j})
+\frac{(1-2y)^2}{4y(1-y)}H(y)\sum\limits_{j=1}^{n-1}\dot\pi(I_{n,j})^2
=\lambda
H(y),
\end{gather*}
for certain representations $\pi\in\hat{\mathrm{SO}}(n)$, including those that are fundamental.

The obvious place to start to look for irreducible representations of ${\mathrm{SO}}(n)$ is among the exterior powers of
the standard representation of ${\mathrm{SO}}(n)$.
It is known that $\Lambda^p(\mathbb{C}^{2\ell})$ are irreducible ${\mathrm{SO}}(2\ell)$-modules for $p=1,\dots,\ell-1$,
and that $\Lambda^\ell(\mathbb{C}^{2\ell})$ splits into the direct sum of two irreducible submodules.
While in the odd case $\Lambda^p(\mathbb{C}^{2\ell+1})$ are irreducible ${\mathrm{SO}}(2\ell+1)$-modules for
$p=1,\dots,\ell$.
See Theorems 19.2 and 19.14 in~\cite{FH91}.

Moreover, $\Lambda^p(\mathbb{C}^{n})$ and $\Lambda^{n-p}(\mathbb{C}^{n})$ are isomorphic ${\mathrm{SO}}(n)$-modules.
In fact, if $\{{\mathbf{e}}_1,\dots,{\mathbf{e}}_n\}$ is the canonical basis of $\mathbb{C}^n$, then the linear map
$\xi:\Lambda^p(\mathbb{C}^{n})\rightarrow\Lambda^{n-p}(\mathbb{C}^{n})$ def\/ined~by
\begin{gather*}
\xi({\mathbf{e}}_{u_1}\wedge\dots\wedge{\mathbf{e}}_{u_p})=(-1)^{u_1+\dots+
u_p}{\mathbf{e}}_{v_1}\wedge\dots\wedge{\mathbf{e}}_{v_{n-p}},
\end{gather*}
where $u_1<\dots<u_p$ and $v_1<\dots<v_{n-p}$ are complementary ordered set of indices, is an
${\mathrm{SO}}(n)$-isomorphism.

All these statements can be established directly upon observing that the elements $I_{ki}=E_{ki}-E_{ik}$ with $1\le
i<k\le n$ form a~basis of the Lie algebra $\mathfrak{so}(n)$, and that
\begin{gather*}
I_{ki}{\mathbf{e}}_k={\mathbf{e}}_i,
\qquad
I_{ki}{\mathbf{e}}_i=-{\mathbf{e}}_k
\qquad
\text{and}
\qquad
I_{ki}{\mathbf{e}}_j=0
\qquad
\text{if}
\quad
j\ne k,i.
\end{gather*}

We will refer to the irreducible ${\mathrm{SO}}(2\ell)$-modules $\Lambda^p(\mathbb{C}^{2\ell})$ for $p=1,\dots,\ell-1$,
respectively, the irreducible ${\mathrm{SO}}(2\ell+1)$-modules $\Lambda^p(\mathbb{C}^{2\ell+1})$ for $p=1,\dots,\ell$,
as the fundamental ${\mathrm{SO}}(2\ell)$-modules, respectively, as the fundamental ${\mathrm{SO}}(2\ell+1)$-modules,
for reasons that will be clarif\/ied in the following Sections~\ref{even} and~\ref{odd}.

\subsection[The even case: $K={\mathrm{SO}}(2\ell)$]{The even case: $\boldsymbol{K={\mathrm{SO}}(2\ell)}$}
\label{even}

First we will study the case $n=2\ell$, with $\ell>2$.
The fundamental weights of $\mathfrak{so}(2\ell,\mathbb{C})$ are
\begin{gather*}
\lambda_p =\epsilon_1+\dots+\epsilon_p,
\qquad
1\le p\le\ell-2,
\\
\lambda_{\ell-1} =\tfrac12(\epsilon_1+\dots+\epsilon_{\ell-1}-\epsilon_\ell),
\qquad
\lambda_{\ell}=\tfrac12(\epsilon_1+\dots+\epsilon_{\ell-1}+\epsilon_\ell).
\end{gather*}

Here we will consider the fundamental $K$-modules
\begin{gather*}
\Lambda^1\big(\mathbb{C}^n\big),\ \Lambda^2\big(\mathbb{C}^n\big),\ \dots, \ \Lambda^{\ell-1}\big(\mathbb{C}^n\big).
\end{gather*}
We will show that the highest weight of $\Lambda^p(\mathbb{C}^n)$ is $\epsilon_1+\dots+\epsilon_p$ for $1\le p\le
\ell-1$.
Observe that $\lambda_{\ell-1}$ and $\lambda_{\ell}$ are not analytically integral and therefore they will not be
considered, although we will also consider the $K$-module with highest weight
$\lambda_{\ell-1}+\lambda_\ell=\epsilon_1+\dots+\epsilon_{\ell-1}$.
Notice that we have already considered the cases $2\lambda_{\ell-1}$ and $2\lambda_{\ell}$ in Section~\ref{mirred},
which are~$M$-irreducible.
We will also show that the fundamental $K$-modules are direct sum of two irreducible~$M$-submodules.

In order to obtain the explicit expression of~$E$ in~\eqref{vectoreq} for a~given irreducible representation~$\pi$ of
$K={\mathrm{SO}}(n)$, of highest weight $\varepsilon_1+\dots+\varepsilon_p$, we are interested to compute
\begin{gather*}
\sum\limits_{j=1}^{n-1}\dot\pi(I_{nj})P_{s}\dot\pi(I_{nj}){\big|_{V_{r}}}=\lambda(r,s)I_{V_{r}},
\end{gather*}
with $r,s=0,1$ corresponding to the two~$M$-submodules $V_0 $ and $V_1$ of the representation~$\pi$, associated with
${\mathbf{m}}_{n-1}=(1,\dots,1,0,\dots,0)\in\mathbb{C}^{\ell-1}$ with $p-1$ and~$p$ ones, respectively (see the
betweenness conditions in Section~\ref{GTB}); being $P_0$ and $P_1$ the respective projections.

Let us consider the standard action of $K={\mathrm{SO}}(n)$ on $V=\mathbb{C}^n$, and take the canonical basis
$\{{\mathbf{e}}_1,\dots,{\mathbf{e}}_n\}$.
Then we have the irreducible $K$-module $\Lambda^p(V)$ for $1\le p\le \ell-1$.
The vector
$({\mathbf{e}}_1-i{\mathbf{e}}_2)\wedge({\mathbf{e}}_3-i{\mathbf{e}}_4)\wedge\dots\wedge({\mathbf{e}}_{2p-1}-i{\mathbf{e}}_{2p})$
is the unique, up to a~scalar, dominant vector and its weight is $(1,\dots,1,0,\dots,0)\in\mathbb{C}^\ell$ with~$p$
ones.
Then, if $V'$ is the subspace generated by $\{{\mathbf{e}}_1,\dots,{\mathbf{e}}_{n-1}\}$, $\Lambda^p(V)$ is the direct
sum of two~$M$-submodules, namely
\begin{gather}
\label{dimpar}
\Lambda^p(V)=V_0\oplus V_1=\Lambda^{p-1}(V')\wedge{\mathbf{e}}_n\oplus \Lambda^p(V')
\end{gather}
whose highest weights are $(1,\dots,1,0,\dots,0)\in\mathbb{C}^{\ell-1}$ with $p-1$ ones and
$(1,\dots,1,0,\dots,0)\in\mathbb{C}^{\ell-1}$ with~$p$ ones, respectively.
It is easy to see that
$({\mathbf{e}}_1-i{\mathbf{e}}_2)\wedge({\mathbf{e}}_3-i{\mathbf{e}}_4)\wedge\dots\wedge({\mathbf{e}}_{2p-3}-i{\mathbf{e}}_{2p-2})\wedge{\mathbf{e}}_n$
is an~$M$-highest weight vector in $\Lambda^{p-1}(V')\wedge{\mathbf{e}}_n$ and that
$({\mathbf{e}}_1-i{\mathbf{e}}_2)\wedge({\mathbf{e}}_3-i{\mathbf{e}}_4)\wedge\dots\wedge({\mathbf{e}}_{2p-1}-i{\mathbf{e}}_{2p})$
is an~$M$ highest weight vector in $\Lambda^p(V')$.

To get $\lambda(0,0)$ it is enough to compute
\begin{gather*}
\sum\limits_{j=1}^{n-1}\dot\pi(I_{nj})P_{0}\dot\pi(I_{nj})({\mathbf{e}}_1\wedge\dots\wedge{\mathbf{e}}_{p-1}\wedge{\mathbf{e}}_n).
\end{gather*}
Since we have that
$\dot\pi(I_{nj})({\mathbf{e}}_1\wedge\dots\wedge{\mathbf{e}}_{p-1}\wedge{\mathbf{e}}_n)
={\mathbf{e}}_1\wedge\dots\wedge{\mathbf{e}}_{p-1}\wedge{\mathbf{e}}_j$
we obtain $P_{0}\dot\pi(I_{nj})({\mathbf{e}}_1\wedge\dots\wedge{\mathbf{e}}_{p-1}\wedge{\mathbf{e}}_n)=0$ and
$\lambda(0,0)=0$.

To get $\lambda(0,1)$ it is enough to compute
\begin{gather*}
\sum\limits_{j=1}^{n-1}\dot\pi(I_{nj})P_{1}\dot\pi(I_{nj})({\mathbf{e}}_1\wedge\dots\wedge{\mathbf{e}}_{p-1}\wedge{\mathbf{e}}_n).
\end{gather*}
We have
\begin{gather*}
P_1\dot\pi(I_{nj})({\mathbf{e}}_1\wedge\dots\wedge{\mathbf{e}}_{p-1}\wedge{\mathbf{e}}_n)=
\begin{cases}
0&
\text{if}\quad
1\le j\le p-1,
\\
{\mathbf{e}}_1\wedge\dots\wedge{\mathbf{e}}_{p-1}\wedge{\mathbf{e}}_j &
\text{if}\quad
p\le j\le n-1.
\end{cases}
\end{gather*}
Therefore we have
\begin{gather*}
\dot\pi(I_{nj})P_1\dot\pi(I_{nj})({\mathbf{e}}_1\wedge\dots\wedge{\mathbf{e}}_{p-1}\wedge{\mathbf{e}}_n)= \begin{cases}
0&\text{if}\quad 1\le j\le p-1,
\\
-{\mathbf{e}}_1\wedge\dots\wedge{\mathbf{e}}_{p-1}\wedge{\mathbf{e}}_n&\text{if}\quad p\le j\le n-1.
\end{cases}
\end{gather*}
Hence $\lambda(0,1)=-(n-p)$.

 Similarly, to get $\lambda(1,0)$ it is enough to compute
\begin{gather*}
\sum\limits_{j=1}^{n-1}\dot\pi(I_{nj})P_{0}\dot\pi(I_{nj})({\mathbf{e}}_1\wedge\dots\wedge{\mathbf{e}}_{p}).
\end{gather*}
We have
\begin{gather*}
\dot\pi(I_{nj})({\mathbf{e}}_1\wedge\dots\wedge{\mathbf{e}}_{p})=
\begin{cases}
-{\mathbf{e}}_1\wedge\dots\wedge{\mathbf{e}}_n\wedge\dots\wedge{\mathbf{e}}_{p}&
\text{if}\quad
1\le j\le p,
\\
0 &
\text{if}\quad
p+1\le j\le n-1,
\end{cases}
\end{gather*}
where ${\mathbf{e}}_n$ appears in the~$j$-place.
Therefore
\begin{gather*}
\dot\pi(I_{nj})P_0\dot\pi(I_{nj})({\mathbf{e}}_1\wedge\dots\wedge{\mathbf{e}}_{p})=
\begin{cases}
-{\mathbf{e}}_1\wedge\dots\wedge{\mathbf{e}}_{p}&\text{if}\quad 1\le j\le p,
\\
0 &\text{if}\quad p+1\le j\le n-1.
\end{cases}
\end{gather*}
Hence $\lambda(1,0)=-p$.

Also it is clear now that
$\sum\limits_{j=1}^{n-1}\dot\pi(I_{nj})P_{1}\dot\pi(I_{nj})({\mathbf{e}}_1\wedge\dots\wedge{\mathbf{e}}_{p})=0$, hence
$\lambda(1,1)=0$.

Therefore, when~$\pi$ is the standard representation of~$K$ in $\Lambda^p(V)$, $1\le p\le \ell-1$, we have
\begin{gather*}
(\lambda(r,s))_{0\le r,s\le1}=
\begin{pmatrix}
0&p-n
\\
-p&0
\end{pmatrix}
.
\end{gather*}

Therefore, we obtain a~more explicit version of Corollary~\ref{eigenfunction2} using~\eqref{vectoreq} and
Remark~\ref{ene}.

\begin{Corollary}
\label{operator2l}
Let~$\Phi$ be an irreducible spherical function on~$G$ of type $\pi\in \hat {\mathrm{SO}}(n)$, $n=2\ell$.
If the highest weight of~$\pi$ is of the form $(1,\dots,1,0,\dots,0)\in\mathbb{C}^\ell$, with~$p$ ones, $1\le p\le
\ell-1$, then the function $H:(0,1)\to{\operatorname{End}}(\mathbb{C}^2)$ associated with~$\Phi$ satisfies
\begin{gather*}
y(1-y)H''(y)+\frac{1}{2}n(1-2y)H'(y)+\frac{1+(1-2y)^2}{4y(1-y)}\left(\begin{matrix}
p-n&0
\\
0&-p
\end{matrix}
\right)H(y)
\\
\qquad
{}+\frac{(1-2y)}{2y(1-y)}\left(\begin{matrix}
0&p-n
\\
-p&0
\end{matrix}
\right)H(y)=\lambda H(y),
\end{gather*}
for some $\lambda\in\mathbb{C}$.
\end{Corollary}

\subsection[The odd case: $K={\mathrm{SO}}(2\ell+1)$]{The odd case: $\boldsymbol{K={\mathrm{SO}}(2\ell+1)}$}
\label{odd}

We now study the case $n=2\ell+1$, with $\ell\ge1$.
The fundamental weights of $\mathfrak{so}(2\ell+1,\mathbb{C})$ are
\begin{gather*}
\lambda_p =\epsilon_1+\dots+\epsilon_p,
\qquad
1\le p\le\ell-1,
\\
\lambda_\ell =\tfrac12(\epsilon_1+\dots+\epsilon_\ell).
\end{gather*}

Here we will consider the fundamental $K$-modules
\begin{gather*}
\Lambda^1\big(\mathbb{C}^n\big), \ \Lambda^2\big(\mathbb{C}^n\big), \ \dots, \ \Lambda^{\ell}\big(\mathbb{C}^n\big).
\end{gather*}
We will show that the highest weight of $\Lambda^p(\mathbb{C}^n)$ is $\epsilon_1+\dots+\epsilon_p$ for $1\le p\le \ell$.
Also we will establish that $\Lambda^p(\mathbb{C}^n)$ splits into the direct sum of two~$M$-submodules for $1\le
p\le\ell-1$, while $\Lambda^\ell(\mathbb{C}^n)$ splits into the sum of three~$M$-submodules; for this reason it will be
treated separately in Section~\ref{mvop}.

Observe that $\lambda_{\ell}$ is not analytically integral and therefore it will not be considered, although we will
consider the $K$-module with highest weight $2\lambda_\ell$.

As in the even case we are interested in computing
\begin{gather*}
\sum\limits_{j=1}^{n-1}\dot\pi(I_{nj})P_{s}\dot\pi(I_{nj})\Big|_{V_{r}}=\lambda(r,s)I_{V_{r}},
\end{gather*}
with $r,s=0,1$ corresponding to the two~$M$-submodules $V_0 $ and $V_1$ of the representation~$\pi$, corresponding to
${\mathbf{m}}_{n-1}=(1,\dots,1,0,\dots,0)\in\mathbb{C}^\ell$ with $p-1$ and~$p$ ones respectively (see the betweenness
conditions in Section~\ref{GTB}).
Being $P_0$ and $P_1$ the respective projections.

Let us consider the standard action of $K={\mathrm{SO}}(n)$ on $V=\mathbb{C}^n$, and take the canonical basis
$\{{\mathbf{e}}_1,\dots,{\mathbf{e}}_n\}$.
Then we have the irreducible $K$-module $\Lambda^p(V)$ for $1\le p\le \ell-1$.
The vector
$({\mathbf{e}}_1-i{\mathbf{e}}_2)\wedge({\mathbf{e}}_3-i{\mathbf{e}}_4)\wedge\dots\wedge({\mathbf{e}}_{2p-1}-i{\mathbf{e}}_{2p})$
is the unique, up to a~scalar, dominant vector and its weight is $(1,\dots,1,0,\dots,0)\in\mathbb{C}^\ell$ with~$p$
ones.
Then, if $V'$ is the subspace generated by $\{{\mathbf{e}}_1,\dots,{\mathbf{e}}_{n-1}\}$, $\Lambda^p(V)$ is the direct
sum of two irreducible~$M$-submodules, namely
\begin{gather}
\label{dimimpar}
\Lambda^p(V)=V_0\oplus V_1=\Lambda^{p-1}(V')\wedge{\mathbf{e}}_n\oplus \Lambda^p(V')
\end{gather}
of highest weights $(1,\dots,1,0,\dots,0)\in\mathbb{C}^{\ell}$ with $p-1$ ones, and
$(1,\dots,1,0,\dots,0)\in\mathbb{C}^{\ell}$ with~$p$ ones, respectively.
It is easy to see that
$({\mathbf{e}}_1-i{\mathbf{e}}_2)\wedge({\mathbf{e}}_3-i{\mathbf{e}}_4)\wedge\dots\wedge({\mathbf{e}}_{2p-3}-i{\mathbf{e}}_{2p-2})\wedge{\mathbf{e}}_n$
is an~$M$-highest weight vector in $\Lambda^{p-1}(V')\wedge{\mathbf{e}}_n$ and that
$({\mathbf{e}}_1-i{\mathbf{e}}_2)\wedge({\mathbf{e}}_3-i{\mathbf{e}}_4)\wedge\dots\wedge({\mathbf{e}}_{2p-1}-i{\mathbf{e}}_{2p})$
is an~$M$ highest weight vector in $\Lambda^p(V')$.

To get $\lambda(0,0)$ it is enough to compute
\begin{gather*}
\sum\limits_{j=1}^{n-1}\dot\pi(I_{nj})P_{0}\dot\pi(I_{nj})({\mathbf{e}}_1\wedge\dots\wedge{\mathbf{e}}_{p-1}\wedge{\mathbf{e}}_n).
\end{gather*}
Since we have that
$\dot\pi(I_{nj})({\mathbf{e}}_1\wedge\dots\wedge{\mathbf{e}}_{p-1}\wedge{\mathbf{e}}_n)
={\mathbf{e}}_1\wedge\dots\wedge{\mathbf{e}}_{p-1}\wedge{\mathbf{e}}_j$,
we obtain $P_{0}\dot\pi(I_{nj})({\mathbf{e}}_1\wedge\dots\wedge{\mathbf{e}}_{p-1}\wedge{\mathbf{e}}_n)=0$ and
$\lambda(0,0)=0$.

To get $\lambda(0,1)$ it is enough to compute
\begin{gather*}
\sum\limits_{j=1}^{n-1}\dot\pi(I_{nj})P_{1}\dot\pi(I_{nj})({\mathbf{e}}_1\wedge\dots\wedge{\mathbf{e}}_{p-1}\wedge{\mathbf{e}}_n).
\end{gather*}
We have
\begin{gather*}
P_1\dot\pi(I_{nj})({\mathbf{e}}_1\wedge\dots\wedge{\mathbf{e}}_{p-1}\wedge{\mathbf{e}}_n)=
\begin{cases}
0&
\text{if}\quad
1\le j\le p-1,
\\
{\mathbf{e}}_1\wedge\dots\wedge{\mathbf{e}}_{p-1}\wedge{\mathbf{e}}_j &
\text{if}\quad
p\le j\le n-1.
\end{cases}
\end{gather*}
Therefore
\begin{gather*}
\dot\pi(I_{nj})P_1\dot\pi(I_{nj})({\mathbf{e}}_1\wedge\dots\wedge{\mathbf{e}}_{p-1}\wedge{\mathbf{e}}_n)=\begin{cases}
0&\text{if}\quad
1\le j\le p-1,
\\
-{\mathbf{e}}_1\wedge\dots\wedge{\mathbf{e}}_{p-1}\wedge{\mathbf{e}}_n&\text{if}\quad
p\le j\le n-1.
\end{cases}
\end{gather*}
Hence $\lambda(0,1)=-(n-p)$.

Similarly, to get $\lambda(1,0)$ it is enough to compute
\begin{gather*}
\sum\limits_{j=1}^{n-1}\dot\pi(I_{nj})P_{0}\dot\pi(I_{nj})({\mathbf{e}}_1\wedge\dots\wedge{\mathbf{e}}_{p}).
\end{gather*}
We have that
\begin{gather*}
\dot\pi(I_{nj})({\mathbf{e}}_1\wedge\dots\wedge{\mathbf{e}}_{p})=
\begin{cases}
-{\mathbf{e}}_1\wedge\dots\wedge{\mathbf{e}}_n\wedge\dots\wedge{\mathbf{e}}_{p}&
\text{if}\quad
1\le j\le p,
\\
0 &
\text{if}\quad
p+1\le j\le n-1,
\end{cases}
\end{gather*}
where ${\mathbf{e}}_n$ appears in the~$j$-place.
Therefore
\begin{gather*}
\dot\pi(I_{nj})P_0\dot\pi(I_{nj})({\mathbf{e}}_1\wedge\dots\wedge{\mathbf{e}}_{p})=
\begin{cases}
-{\mathbf{e}}_1\wedge\dots\wedge{\mathbf{e}}_{p} &\text{if}\quad
1\le j\le p,
\\
0 &\text{if}\quad
p+1\le j\le n-1.
\end{cases}
\end{gather*}
Hence $\lambda(1,0)=-p$.

Also it is clear now that
$\sum\limits_{j=1}^{n-1}\dot\pi(I_{nj})P_{1}\dot\pi(I_{nj})({\mathbf{e}}_1\wedge\dots\wedge{\mathbf{e}}_{p})=0$, hence
$\lambda(1,1)=0$.

Therefore, when~$\pi$ is the standard representation of~$K$ in $\Lambda^p(V)$, $1\le p\le \ell-1$, we have
\begin{gather*}
(\lambda(r,s))_{0\le r,s\le1}=
\begin{pmatrix}
0&p-n
\\
-p&0
\end{pmatrix}
.
\end{gather*}

Therefore, we obtain a~more explicit version of Corollary~\ref{eigenfunction2} using~\eqref{vectoreq} and
Remark~\ref{ene}.
\begin{Corollary}
\label{operator2l+1}
Let~$\Phi$ be an irreducible spherical function on~$G$ of type $\pi\in \hat {\mathrm{SO}}(n)$, $n=2\ell+1$.
If the highest weight of~$\pi$ is of the form $(1,\dots,1,0,\dots,0)$ $\in\mathbb{C}^\ell$, with~$p$ ones, $1\le p\le
\ell-1$, then the function $H:(0,1)\to{\operatorname{End}}(\mathbb{C}^2)$ associated with~$\Phi$ satisfies
\begin{gather*}
y(1-y)H''(y)+\frac{1}{2}n(1-2y)H'(y)+\frac{1+(1-2y)^2}{4y(1-y)}\left(\begin{matrix}
p-n&0
\\
0&-p
\end{matrix}
\right)H(y)
\\
\qquad
{}+\frac{(1-2y)}{2y(1-y)}\left(\begin{matrix}
0&p-n
\\
-p&0
\end{matrix}
\right)H(y)=\lambda H(y),
\end{gather*}
for some $\lambda\in\mathbb{C}$.
\end{Corollary}

\section[The spherical functions of fundamental $K$-types]{The spherical functions of fundamental $\boldsymbol{K}$-types}
\label{s2l}

Let $n=2\ell$, the irreducible spherical functions of $K$-type
\begin{gather*}
{\mathbf{m}}_{n}=(1,\dots,1,0,\dots,0)\in\mathbb{C}^\ell,
\end{gather*}
with~$p$ ones, $1\le p\le\ell-1$, are those associated with the irreducible representations of~$G$ of highest weights of
the form ${\mathbf{m}}_{n+1}=(w+1,1,\dots,1,\delta,0,\dots,0)$ $\in\mathbb{C}^\ell$ that interlaces ${\mathbf{m}}_n$,
\[
\begin{array}{@{}cccccccccccccccccc} w+1 &{} &1 &\dots &1 &{} &\delta &{} &0 &\dots &0 &{}
\\
{} &1 &\dots &{} &\dots &1 &{} &0 &\dots &{} &\dots &0
\end{array}.
\]

We now consider the $K$-module $\Lambda^p(\mathbb{C}^{n})$ which has highest weight ${\mathbf{m}}_{n}$.

For $w=0$ and $\delta=0$ we consider the~$G$-module $\Lambda^p(\mathbb{C}^{n+1})$ whose highest weight is
${\mathbf{m}}_{n+1}$, and we have the following $K$-module decomposition
\begin{gather*}
\Lambda^p\big(\mathbb{C}^{n+1}\big) =\Lambda^p(\mathbb{C}^n)\oplus \Lambda^{p-1}\big(\mathbb{C}^n\big)\wedge{\mathbf{e}}_{n+1},
\end{gather*}
where $\Lambda^p(\mathbb{C}^{n})$ is the sum of two ${\mathrm{SO}}(n-1)$-modules:
\begin{gather*}
\Lambda^p\big(\mathbb{C}^{n}\big)=\Lambda^p\big(\mathbb{C}^{n-1}\big)\oplus\Lambda^{p-1}\big(\mathbb{C}^{n-1}\big)\wedge{\mathbf{e}}_n.
\end{gather*}
We observe that
\begin{gather*}
a(s)({\mathbf{e}}_1\wedge\dots\wedge{\mathbf{e}}_{p-1}\wedge{\mathbf{e}}_n)={\mathbf{e}}_1\wedge\dots\wedge{\mathbf{e}}_{p-1}\wedge(\cos
s {\mathbf{e}}_n-\sin s{\mathbf{e}}_{n+1})
\\
\qquad
=\cos s({\mathbf{e}}_1\wedge\dots\wedge{\mathbf{e}}_{p-1}\wedge{\mathbf{e}}_n)-\sin
s({\mathbf{e}}_1\wedge\dots\wedge{\mathbf{e}}_{p-1}\wedge{\mathbf{e}}_{n+1}).
\end{gather*}
Hence, if $\Phi_{0}$ is the spherical function associated with the irreducible representation of~$G$ of highest weight
${\mathbf{m}}_{n+1}=(1,1,\dots,1,\delta,0,\dots,0)\in\mathbb{C}^\ell$ with $\delta=0$, we have that
\begin{gather*}
\Phi_0(a(s))({\mathbf{e}}_1\wedge\dots\wedge{\mathbf{e}}_{p-1}\wedge{\mathbf{e}}_n)=\cos
s({\mathbf{e}}_1\wedge\dots\wedge{\mathbf{e}}_{p-1}\wedge{\mathbf{e}}_n).
\end{gather*}
Also we have that $a(s)({\mathbf{e}}_1\wedge\dots\wedge{\mathbf{e}}_p)={\mathbf{e}}_1\wedge\dots\wedge{\mathbf{e}}_p$.
Thus the vector valued function $F_0(s)$ given by the irreducible spherical function $\Phi_0(a(s))$ is
\begin{gather*}
F_0(s)=\left(
\begin{matrix}
\cos s
\\
1
\end{matrix}
\right).
\end{gather*}

For $w=0$ and $\delta=1$ we consider the~$G$-module $\Lambda^{p+1}(\mathbb{C}^{n+1})$ whose highest weight
${\mathbf{m}}_{n+1}$, and for $1\le p\le\ell-1$ we have the following $K$-module decomposition
\begin{gather*}
\Lambda^{p+1}\big(\mathbb{C}^{n+1}\big) =\Lambda^{p+1}\big(\mathbb{C}^n\big)\oplus \Lambda^{p}\big(\mathbb{C}^n\big)\wedge{\mathbf{e}}_{n+1},
\end{gather*}
where $\Lambda^{p}(\mathbb{C}^n)\wedge{\mathbf{e}}_{n+1}$ is the sum of two ${\mathrm{SO}}(n-1)$-modules:
\begin{gather*}
\Lambda^{p}\big(\mathbb{C}^n\big)\wedge{\mathbf{e}}_{n+1}=\Lambda^{p}\big(\mathbb{C}^{n-1}\big)\wedge{\mathbf{e}}_{n+1}
\oplus\Lambda^{p-1}\big(\mathbb{C}^{n-1}\big)\wedge{\mathbf{e}}_n\wedge{\mathbf{e}}_{n+1}.
\end{gather*}
We observe that
\begin{gather*}
a(s)({\mathbf{e}}_1\wedge\dots\wedge{\mathbf{e}}_{p-1}\wedge{\mathbf{e}}_n\wedge{\mathbf{e}}_{n+1})=
{\mathbf{e}}_1\wedge\dots\wedge{\mathbf{e}}_{p-1}\wedge(\sin s {\mathbf{e}}_n+\cos s{\mathbf{e}}_{n+1})
\\
\qquad
=\sin s({\mathbf{e}}_1\wedge\dots\wedge{\mathbf{e}}_{p-1}\wedge{\mathbf{e}}_n)+\cos s
({\mathbf{e}}_1\wedge\dots\wedge{\mathbf{e}}_{p-1}\wedge{\mathbf{e}}_{n+1}).
\end{gather*}

Hence, if $\Phi_{1}$ is the spherical function associated with the irreducible representation of~$G$ of highest weight
${\mathbf{m}}_{n+1}=(1,1,\dots,1,\delta,0,\dots,0)\in\mathbb{C}^\ell$ with $\delta=1$, we have that $
\Phi_1(a(s))({\mathbf{e}}_1\wedge\dots\wedge{\mathbf{e}}_{p-1}\wedge{\mathbf{e}}_n\wedge{\mathbf{e}}_{n+1})=\cos
s({\mathbf{e}}_1\wedge\dots\wedge{\mathbf{e}}_{p-1}\wedge{\mathbf{e}}_n\wedge{\mathbf{e}}_{n+1})$.
Also we have that
\begin{gather*}
a(s)({\mathbf{e}}_1\wedge\dots\wedge{\mathbf{e}}_p\wedge{\mathbf{e}}_{n+1})
={\mathbf{e}}_1\wedge\dots\wedge{\mathbf{e}}_{p-1}\wedge{\mathbf{e}}_n\wedge{\mathbf{e}}_{n+1}.
\end{gather*}

Thus the vector valued function $F_1(s)$ given by the irreducible spherical function $\Phi_1(a(s))$ is
\begin{gather*}
F_1(s)=\left(
\begin{matrix}
1
\\
\cos s
\end{matrix}
\right).
\end{gather*}

\begin{Definition}
We shall consider the $2\times2$ matrix-valued function $\Psi=\Psi(y)$, for $0<y<1$, whose columns are given by the
functions $H_0(y)=F_0(s)$ and $H_1(y)=F_1(s)$, with $\cos s=2y-1$:
\begin{gather}
\label{Psi}
\Psi(y)=\left(
\begin{matrix}
2y-1&1
\\
1&2y-1
\end{matrix}
\right).
\end{gather}
\end{Definition}

Since the functions $H_0(y)$ and $H_1(y)$ are associated with irreducible spherical functions, they satisfy the
dif\/ferential equation given in Corollary~\ref{operator2l}; moreover, the respective eigenvalues are $\lambda=-p$ and
$\lambda=p-n$.
Therefore, we have
\begin{gather*}
y(1-y)\Psi''+\frac{1}{2}n(1-2y)\Psi'+\frac{1+(1-2y)^2}{4y(1-y)}\left(
\begin{matrix}
p-n&0
\\
0&-p
\end{matrix}
\right)\Psi
\\
\qquad
{}+\frac{(1-2y)}{2y(1-y)}\left(
\begin{matrix}
0&p-n
\\
-p&0
\end{matrix}
\right)\Psi =\Psi\left(
\begin{matrix}
-p&0
\\
0&p-n
\end{matrix}
\right).
\end{gather*}

Furthermore, it is easy to check that the function $\Psi(y)$ also satisfy the equation above even when~$n$ is odd.

\begin{Theorem}
\label{ecpar}
The function~$\Psi$ can be used to obtain a~hypergeometric differential equation from the one given in
Corollaries~{\rm \ref{operator2l}} and~{\rm \ref{operator2l+1}}.
Precisely, if~$H$ is a~vector-valued solution of the differential equation in Corollaries~{\rm \ref{operator2l}}
or~{\rm \ref{operator2l+1}}, with eigenvalue~$\lambda$, then $P=\Psi^{-1}H$ is a~solution of $DP=\lambda P$, where~$D$ is the
hypergeometric differential operator given~by
\begin{gather*}
DP=y(1-y)P''-\left(
\begin{matrix}
(\frac n2+1)(2y-1)&-1
\\
-1&(\frac n2+1)(2y-1)
\end{matrix}
\right)P' -\left(
\begin{matrix}
p&0
\\
0&n-p
\end{matrix}
\right)P.
\end{gather*}
\end{Theorem}

\begin{proof}
By hypothesis we have that
\begin{gather*}
y(1-y)H''(y)+\frac{1}{2}n(1-2y)H'(y)+\frac{1+(1-2y)^2}{4y(1-y)}\left(\begin{matrix}
p-n&0
\\
0&-p
\end{matrix}
\right)H(y)
\\
\qquad
{}+\frac{(1-2y)}{2y(1-y)}\left(\begin{matrix}
0&p-n
\\
-p&0
\end{matrix}
\right)H(y)=\lambda H(y),
\end{gather*}
Then, writing $H=\Psi P$, we have
\begin{gather*}
y(1-y)P''+\big(2y(1-y)\Psi^{-1}\Psi'+\frac n2(1-2y)I\big)P'
\\
\qquad
{}+\Psi^{-1}\left(y(1-y)\Psi''+\frac n2(1-2y)\Psi'+\frac{1+(1-2y)^2}{4y(1-y)}\left(
\begin{matrix}
p-n&0
\\
0&-p
\end{matrix}
\right)\Psi\right.
\\
\qquad
{}+\left.\frac{(1-2y)}{2y(1-y)}
\begin{pmatrix}
0&p-n
\\
-p&0
\end{pmatrix}
\Psi\right)P=\lambda P.
\end{gather*}
Now we compute
\begin{gather*}
2y(1-y)\Psi^{-1}\Psi'=\frac{4y(1-y)}{4y(y-1)}\left(
\begin{matrix}
2y-1&-1
\\
-1&2y-1
\end{matrix}
\right)=-\left(
\begin{matrix}
2y-1&-1
\\
-1&2y-1
\end{matrix}
\right).
\end{gather*}
Therefore
\begin{gather*}
y(1-y)P''-\left(
\begin{matrix}
(\frac n2+1)(2y-1)&-1
\\
-1&(\frac n2+1)(2y-1)
\end{matrix}
\right)P' -\left(
\begin{matrix}
\lambda+p&0
\\
0&\lambda+n-p
\end{matrix}
\right)P=0.
\end{gather*}
This completes the proof of the theorem.
\end{proof}

\subsection[$\Delta$-eigenvalues of spherical functions]{$\boldsymbol{\Delta}$-eigenvalues of spherical functions}

As we said, when $n=2\ell$ the irreducible spherical functions of the pair $({\mathrm{SO}}(n+1),{\mathrm{SO}}(n))$, of
type ${\mathbf{m}}_n=(1,\dots,1,0\dots,0)\in\mathbb{C}^\ell$ with~$p$ ones, $1\leq p\leq\ell-1$ are those associated
with the irreducible representations~$\tau$ of~$G$ of highest weights of the form
${\mathbf{m}}_{n+1}=(w+1,1,\dots,1,\delta,0,\dots,0)\in\mathbb{C}^\ell$ with $p-1$ ones, such that the following pattern
holds
\begin{gather*}
\begin{array}{@{}cccccccccccccccccccc}
w+1 &{} &1 &\dots &1 &{} &\delta &{} &0 &\dots &0&{}
\\
{} &1 &\dots &{} &\dots &1 &{} &0 &\dots &{} &\dots &0
\end{array}.
\end{gather*}

Let $\Phi_{w,\delta}$ be the corresponding spherical function.
Then $\Delta \Phi_{w,\delta}=\lambda\Phi_{w,\delta}$, where the eigenvalue $\lambda=\lambda_n(w,\delta)$ can be computed
from the expression $\Delta=Q_{n+1}-Q_n$.
If $v\in V_{{\mathbf{m}}_{n+1}}$ is a~highest weight vector from~\eqref{qvimpar} we have
\begin{gather*}
\dot\tau(Q_{2\ell+1})v =-\big((w+1)^2+(2\ell-1)(w+1)+(2\ell-p)(p-1)+2\delta(\ell-p)\big)v.
\end{gather*}
If $v\in V_{{\mathbf{m}}_{2\ell}}$ is a~highest weight vector, then from~\eqref{qvpar} we have
\begin{gather*}
\dot\pi(Q_n)v=-p(2\ell-p)v.
\end{gather*}
Since $\Delta=Q_{n+1}-Q_n$ it follows that
\begin{gather*}
\lambda_{2\ell}(w,\delta) =-(w+1)^2-(2\ell-1)(w+1)+(2\ell-p)-2\delta(\ell-p)
\end{gather*}
Analogously, we obtain that the eigenvalues of the spherical functions $\Phi_{w,\delta}$ of the pair
$({\mathrm{SO}}(2\ell+2),{\mathrm{SO}}(2\ell+1))$ are of the form
\begin{gather*}
\lambda_{2\ell+1}(w,\delta) =-(w+1)(w+2\ell+1)+2\ell-p+1-\delta2(\ell-p)-\delta^2,
\end{gather*}
here~$\delta$ is $0$ or $1$ when we are in the cases $1\le p<\ell$ but~$\delta$ could also be $-1$ in the particular
case $p=\ell$.

Therefore, we have that the eigenvalues of the spherical functions $\Phi_{w,\delta}$ of the pair
$({\mathrm{SO}}(n+1),{\mathrm{SO}}(n))$ are of the form
\begin{gather}
\label{lambda}
\lambda_{n}(w,\delta)=
\begin{cases}
-w(w+n+1)-p
&\text{if}\quad
\delta=0,
\\
-w(w+n+1)-n+p
&\text{if}\quad
\delta=\pm1.
\end{cases}
\end{gather}

\subsection[Polynomial eigenfunctions of the hypergeometric operator $D$]{Polynomial eigenfunctions of the hypergeometric operator $\boldsymbol{D}$}

Let $D $ be the dif\/ferential operator on the real line introduced in Theorem~\ref{ecpar}:
\begin{gather}
\label{D}
DP=y(1-y)P''+(C-yU)P'-VP,
\end{gather}
with
\begin{gather*}
C=
\begin{pmatrix}
(n/2+1)&1
\\
1&(n/2+1)
\end{pmatrix},
\qquad
U=(n+2)I,
\qquad
V= \left(
\begin{matrix}
p&0
\\
0&n-p
\end{matrix}
\right),
\end{gather*}
where~$n$ is of the form $2\ell$ or $2\ell+1$ for $\ell\in\mathbb{N}$ and $1\le p<\ell$.

We will study the $\mathbb{C}^{2}$-vector valued polynomial eigenfunctions of~$D$.

The equation $DP=\lambda P$ is an instance of a~matrix hypergeometric dif\/ferential equation studied in~\cite{T03}.
Since the eigenvalues of~$C$, $n/2$ and $n/2+2$, are not in $-\mathbb{N}_0$ the function~$P$ is determined~by $P_0=P(0)$.
For $|y|<1$ it is given~by
\begin{gather*}
P(y)={}_2H_1\left(\begin{matrix}
U,V+\lambda
\\
C
\end{matrix}
; y\right)P_0=\sum\limits_{j=0}^{\infty}\frac{y^j}{j!} [C;U;V+\lambda]_j P_0,
\qquad
P_0\in \mathbb{C}^{2},
\end{gather*}
where the symbol $[C;U;V+\lambda]_j$ is inductively def\/ined~by
\begin{gather*}
[C;U;V+\lambda]_0  =1,
\\
[C;U;V+\lambda]_{j+1}  = \left(C+j\right)^{-1}(j(U+j-1)+V+\lambda)[C;U;V+\lambda]_j,
\end{gather*}
for all $j\geq 0$.

Therefore, we have that there exists a~polynomial solution if and only if the coef\/f\/icient $[C;U;V+\lambda]_{j+1}$ is
a~singular matrix for some $j\in \mathbb{Z}$.
Since the matrix $C+j$ is invertible for all $j\in\mathbb{N}_0$, we have that there is a~polynomial solution of
degree~$j$ for $DP=\lambda P$ if and only if there exists $P_0\in\mathbb{C}^{2}$ such that $[C;U;V+\lambda]_{j}P_0\neq
0$ and $(j(U+j-1)+V+\lambda)[C;U;V+\lambda]_{j}P_0=0$.

Now we easily observe that the only possible values for~$\lambda$ such that $j(U+j-1)+V+\lambda$ has non trivial kernel
are those given in~\eqref{lambda}.
Then, if $\lambda=-w(w+n+1)-p$, it is easy to check that the f\/irst and only~$j$ for which $j(U+j-1)+V+\lambda$ is
singular is $j=w$, and its kernel (of dimension~$1$) is the subspace generated by $\left(\begin{smallmatrix}
1
\\
0
\end{smallmatrix}\right)$.
Analogously, if $\lambda=-w(w+n+1)-n+p$, it is easy to check that the f\/irst and only~$j$ for which $j(U+j-1)+V+\lambda$
is singular is $j=w$, and its kernel (of dimension~$1$) is the subspace generated by $\left(\begin{smallmatrix}
0
\\
1
\end{smallmatrix}\right)$ respectively.
Therefore we have the following result.
\begin{Theorem}
\label{polsol2l}
For a~given $\ell\in\mathbb{N}$ take $n=2\ell$ or $2\ell+1$ and $1\le p\le\ell-1$, then the polynomial eigenfunctions of
\begin{gather*}
DP=y(1-y)P''+(C-yU)P'-VP,
\end{gather*}
with
\begin{gather*}
C=
\begin{pmatrix}
(n/2+1)&1
\\
1&(n/2+1)
\end{pmatrix},
\qquad
U=(n+2)I,
\qquad
V= \left(
\begin{matrix}
p&0
\\
0&n-p
\end{matrix}
\right)
\end{gather*}
have eigenvalues $-w(w+n+1)-p$ or $-w(w+n+1)-n+p$, with $w\in\mathbb{N}_0$; in both cases the degree of the polynomial is~$w$ with leading coefficient a multiple of $\left(\begin{smallmatrix} 1\\ 0 \end{smallmatrix}\right)$ or $\left(\begin{smallmatrix} 0\\ 1 \end{smallmatrix}\right)$, respectively.
\end{Theorem}

\section{The inner product}
\label{innerprod}

Given a~f\/inite dimensional irreducible representation~$\pi$ of $ K$ in the vector space $V_\pi$ let $(C(G)\otimes
{\operatorname{End}} (V_\pi))^{K\times K}$ be the space of all continuous functions $\Phi: G\longrightarrow
{\operatorname{End}}(V_\pi)$ such that $\Phi(k_1gk_2)=\pi(k_1)\Phi(g)\pi(k_2)$ for all $g\in G$, $k_1,k_2\in K$.
Let us equip $V_\pi$ with an inner product such that $\pi(k)$ becomes unitary for all $k\in K$.
Then we introduce an inner product in the vector space $(C(G)\otimes {\operatorname{End}} (V_\pi))^{K\times K}$~by
def\/ining
\begin{gather*}
\langle \Phi_1,\Phi_2 \rangle =\int_{G} \operatorname{tr} (\Phi_1(g)\Phi_2(g)^*) dg,
\end{gather*}
where $dg$ denote the Haar measure on $ G$ normalized by $\int_G dg=1$, and where $\Phi_2(g)^*$ denotes the adjoint of
$\Phi_2(g)$ with respect to the inner product in $V_\pi$.

By using Schur's orthogonality relations for the unitary irreducible representations of $ G$, it follows that if
$\Phi_1$ and $\Phi_2$ are non equivalent irreducible spherical functions, then they are orthogonal with respect to the
inner product $\langle\cdot,\cdot\rangle$, i.e.~
\begin{gather*}
\langle \Phi_1,\Phi_2 \rangle =0.
\end{gather*}

Recall that, given an irreducible spherical function~$\Phi$ of type~$\pi$ of the pair $(G,K)$, the function $\Phi(a(s))$
is scalar valued when restricted to any ${\mathrm{SO}}(n-1)$-module (see~\eqref{a} for $a(s)$).
We shall denote by~$m$ the number of ${\mathrm{SO}}(n-1)$-submodules of~$\pi$, and by $d_1,d_2,\dots,d_m$ the respective
dimensions of each one of those submodules.

In particular, if $\Phi_1$ and $\Phi_2$ are two irreducible spherical functions of type $\pi\in\hat K$, we consider the
vector valued functions $H_1(y)$ and $H_2(y)$ given by the diagonal matrix valued functions $\Phi_1 (a (s))$ and $\Phi_2
(a (s))$ (see Remark~\ref{hache}), with $y=(\cos s+1)/2$, respectively, denoting
\begin{gather*}
H_1(y)=(h_1(y),\dots, h_m(y))^t,
\qquad
H_2(y)=(f_1(y),\dots, f_m(y))^t.
\end{gather*}

\begin{Proposition}
\label{prodint}
If $\Phi_1$, $\Phi_2$ are two irreducible spherical functions of type $\pi\in\hat K$ then
\begin{gather*}
\langle \Phi_1,\Phi_2\rangle =\frac{(n-1)!!}{(n-2)!!}
\frac{2}{\omega_*} \sum\limits_{i=1}^m d_i \int_{0}^{1}(y(1- y))^{n/2-1}  h_i(y)\overline{f_i(y)}dy,
\end{gather*}
with $\omega_*=\pi$ if~$n$ is even and $\omega_*=2$ if~$n$ is odd.
\end{Proposition}

\begin{proof}
Let $A=\exp \mathbb{R} I_{n+1,n}$ be the Lie subgroup of~$G$ of all elements of the form{\samepage
\begin{gather*}
a(s)=\exp s I_{n+1,n} =
\begin{pmatrix}
I_{n-1}&{\bf 0}&{\bf 0}
\\
{\bf 0} &\cos s&\sin s
\\
{\bf 0} &-\sin s&\cos s
\end{pmatrix}
,
\qquad
s\in \mathbb{R},
\end{gather*}
where $I_{n-1}$ denotes the identity matrix of size $n-1$.}

Now \cite[Theorem 5.10, p.~190]{H00} establishes that for every $f\in C(G/K)$ and a~suitable cons\-tant~$c_*$
\begin{gather*}
\int_{G/K} f(gK)dg_K=c_*\int_{K/M}\left(\int_{-\pi}^{\pi} \delta_*(a(s))f(ka(s)K)ds\right)dk_M,
\end{gather*}
where $dg_K$ and $dk_M$ are respectively the invariant measures on $G/K$ and $K/M$ normalized by $\int_{G/K}
dg_K=\int_{K/M} dk_M=1$ and the function $\delta_*:A\longrightarrow \mathbb{R} $ is def\/ined~by
\begin{gather*}
\delta_*(a(s))=\prod\limits_{\nu\in\Sigma^+} |\sin is \nu(I_{n+1,n})|,
\end{gather*}
with $\Sigma^+$ the set of those positive roots whose restrictions to $\mathfrak{a}$, the Lie algebra of~$A$, are not
zero.
In our case we have $\delta_*(a(s))=|\sin^{n-1}s| $.

To f\/ind the value of $c_*$ we consider the function $f\equiv1$, having then
\begin{gather*}
1=2c_* \int_{0}^{\pi}\sin^{n-1}s ds.
\end{gather*}
Since
\begin{gather*}
\int\sin^{n-1}s ds=-\frac{1}{n-1} \sin^{n-2}s \cos s+\frac{n-2}{n-1}\int \sin^{n-3} ds,
\end{gather*}
we obtain that, for $n=2\ell$ or $2\ell+1$,
\begin{gather*}
\int_0^\pi\sin^{n-1}sds=\frac{n-2}{n-1}\frac{n-4}{n-3}\cdots\frac{n-2\ell+1}{n-2\ell+2}\int_0^\pi\sin^{n-2\ell}s ds.
\end{gather*}
Therefore{\samepage
\begin{gather*}
c_*=\frac{(n-1)!!}{(n-2)!!}
\frac{1}{2\omega_*},
\end{gather*}
with $\omega_*=\pi$ for $n=2\ell$ and $\omega_*=2$ for $2\ell+1$.}

Since the function $g\mapsto \operatorname{tr}(\Phi_1(g)\Phi_2(g)^*)$ is invariant under left and right multiplication
by elements in~$K$, we have
\begin{gather*}
\langle \Phi_1,\Phi_2\rangle =\int_G\operatorname{tr} (\Phi_1(g)\Phi_2(g)^*)dg =2 c_* \int_{0}^{\pi} \sin^{n-1}
s\operatorname{tr}\left(\Phi_1(a(s)\Phi_2(a(s))^*)\right)ds.
\end{gather*}

If we put $y=\tfrac12(\cos s +1)$ for $0<s<\pi$ we have
\begin{gather*}
\operatorname{tr}\left(\Phi_1(a(s)\Phi_2(a(s))^*)\right)
=\sum\limits_{i=1}^m d_ih_i(y)\overline{f_i(y)}.
\end{gather*}

Then
\begin{gather*}
\langle \Phi_1,\Phi_2\rangle =4c_* \sum\limits_{i=1}^m d_i \int_{0}^{1}(4y(1- y))^{(n-2)/2} h_i(y)\overline{f_i(y)}dy,
\end{gather*}
and the proposition follows.
\end{proof}

\begin{Proposition}
\label{Deltasim}
If $\Phi_1,\Phi_2\in(C^\infty(G)\otimes {\operatorname{End}} (V_\pi))^{K\times K}$ then
\begin{gather*}
\langle \Delta\Phi_1, \Phi_2\rangle=\langle \Phi_1,\Delta \Phi_2\rangle.
\end{gather*}
\end{Proposition}
\begin{proof}
If we apply a~left invariant vector f\/ield $X\in \mathfrak g$, to the function on~$G$ given~by
$g\mapsto\operatorname{tr}(\Phi_1(g)\Phi_2(g)^*)$, and then we integrate over~$G$ we obtain
\begin{gather*}
0= \int_G \operatorname{tr}\left((X\Phi_1)(g)\Phi_2(g)^*\right) dg+\int_G
\operatorname{tr}\left(\Phi_1(g)(X\Phi_2)(g)^*\right) dg.
\end{gather*}
Therefore $\langle X\Phi_1,\Phi_2\rangle=-\langle \Phi_1,X\Phi_2\rangle$.
Now let $\tau:\mathfrak g_\mathbb{C}\longrightarrow \mathfrak g_\mathbb{C}$ be the conjugation of $\mathfrak
g_\mathbb{C}$ with respect to the real linear form $\mathfrak g$.
Then $-\tau$ extends to a~unique antilinear involutive ${}^*$ operator on $D(G)$ such that
$(D_1D_2)^*=D_2^*D_1^*$ for all $D_1,D_2\in D(G)$.
This follows easily from the fact that the universal enveloping algebra over $\mathbb{C}$ of $\mathfrak g$ is
canonically isomorphic to $D(G)$.
Then it follows that $\langle D\Phi_1, \Phi_2\rangle=\langle \Phi_1,D^* \Phi_2\rangle$.

Finally, it is easy to verify that $\Delta^*=\Delta$.
\end{proof}

\subsection[Spherical functions as polynomial solutions of $DP=\lambda P $]{Spherical functions as polynomial solutions of $\boldsymbol{DP=\lambda P}$}

Let us consider $\widetilde D $, the dif\/ferential operator on $(0,1)$ introduced in Corollaries~\ref{operator2l}
and~\ref{operator2l+1}:
\begin{gather}
y(1-y)H''(y)+\frac{1}{2}n(1-2y)H'(y)+\frac{1+(1-2y)^2}{4y(1-y)}\left(\begin{matrix}
p-n&0
\\
0&-p
\end{matrix}
\right)H(y)
\nonumber\\
\qquad
{}+\frac{(1-2y)}{2y(1-y)}\left(\begin{matrix}
0&p-n
\\
-p&0
\end{matrix}
\right)H(y) =\lambda H(y),\label{Dtilde}
\end{gather}
Recall that the operator~$D$ that appears in~\eqref{D} extends the dif\/ferential operator $D=\Psi \widetilde D \Psi^{-1}$
to the whole real line, where
\begin{gather*}
\Psi(y)=\left(
\begin{matrix}
2y-1&1
\\
1&2y-1
\end{matrix}
\right)
\end{gather*}
is the matrix function given in~\eqref{Psi} and used in Theorem~\ref{ecpar}.

We want to focus our attention on the following vector spaces of $\mathbb{C}^{2}$-valued analytic functions on $(0,1)$:
\begin{gather*}
S_\lambda= \big\{H=H(y):\widetilde D H=\lambda H,~H(\tfrac{\cos s+1}{2})~\text{analytic at}~s=0\big\},
\\
W_\lambda= \big\{P=P(y):DP=\lambda P,~\text{analytic on}~[0,1]\big\}.
\end{gather*}

From Theorem~\ref{ecpar} we know that the correspondence $P\mapsto\Psi P$ is an injective linear map from~$W_\lambda$
into~$S_\lambda$.
Now we want to prove that this map is bijective.

\begin{Theorem}
\label{isomo}
The linear map $P\mapsto\Psi P$ is an isomorphism from $W_\lambda$ onto $S_\lambda$.
\end{Theorem}
\begin{proof}

A vector valued function $P\in W_\lambda$ is an eigenfunction of the hypergeometric operator~$D$.
Since it is analytic at $y=1$ it is determined by $P(1)$, therefore $\dim(W_\lambda)=2$.

On the other hand, if $H\in S_\lambda$ then there is a~function $F(s)$ analytic at $s=0$, such that it extends the
function $H(\frac{\cos s +1}{2})$ def\/ined on $(0,\pi)$.
Then,~$F$ satisf\/ies the following dif\/ferential equation
\begin{gather*}
F''(s)+(n-1)\frac{\cos s}{\sin s}F'(s)+\frac{1+\cos^2s}{\sin^2s}\left(\begin{matrix}
p-n & 0
\\
0 & -p
\end{matrix}
\right)F(s) \\
\qquad{} -2\frac{\cos s}{\sin^2 s}\left(\begin{matrix}
0&p-n
\\
-p&0
\end{matrix}
\right) F(s)=\lambda F(s),
\end{gather*}
or equivalently
\begin{gather}
\sin^2 s F''(s)+\frac{n-1}{2}\sin (2s)F'(s)+(2-\sin^2s)\left(\begin{matrix}
p-n & 0
\\
0 & -p
\end{matrix}
\right)F(s)
\nonumber\\
\qquad
{}-2\cos s \left(\begin{matrix}
0&p-n
\\
-p&0
\end{matrix}
\right) F(s)=\lambda \sin^2 s F(s),\label{21}
\end{gather}

Let $a_j\in\mathbb{C}^2$ and $\alpha_j,\beta_j, \gamma_j \in \mathbb{C}$, for $j\ge 0$, be the Taylor coef\/f\/icients
of~$F$, $\sin$, $\sin^2$ and $\cos$ at $s=0$:
\begin{gather*}
F(s)=\sum\limits_{j\ge 0} a_j s^j,
\qquad
\sin s=\sum\limits_{j\ge 1} \alpha_j s^j,
\\
F'(s)=\sum\limits_{j\ge 0} a_{j+1}(j+1) s^j,
\qquad
\sin^2 s=\sum\limits_{j\ge 2} \beta_j s^j,
\\
F''(s)=\sum\limits_{j\ge 0} a_{j+2}(j+2)(j+1) s^j,
\qquad
\cos s=\sum\limits_{j\ge 0} \gamma_j s^j.
\end{gather*}

Therefore, from~\eqref{21} we have
\begin{gather*}
\sum\limits_{j\ge 0}\left[\sum\limits_{k=0}^{j-2}\beta_{j-k}
a_{k+2}(k+2)(k+1)+\frac{n-1}{2}\sum\limits_{k=0}^{j-1}2^{j-k}\alpha_{j-k} a_{k+1}(k+1)+\left(\begin{matrix}
p-n&0
\\
0&-p
\end{matrix}
\right) \right.
\\
\qquad
\times \left.\left(2 a_j -\sum\limits_{k=0}^{j-2}\beta_{j-k} a_k \right)-2\left(\begin{matrix}
0&p-n
\\
-p&0
\end{matrix}
\right)\sum\limits_{k=0}^{j}\gamma_{j-k} a_k \right]s^j=\lambda \sum\limits_{j\ge 0}\left[\sum\limits_{k=0}^{j-2}
\beta_{j-k} a_k\right]s^j.
\end{gather*}
Hence, since $\beta_2=\alpha_1=	\gamma_0=1$, we have that
\begin{gather*}
\left[j(j-1) +(n-1) j+ 2 \left(
\begin{matrix}
p-n&-p+n
\\
p&-p
\end{matrix}
\right) \right]a_j
\end{gather*}
is a~linear combination with matrix coef\/f\/icients of $\{a_0, a_1, \dots, a_{j-1} \}$; it is clear that for $j=1$ and
$j>2$ the matrix above is non singular, therefore $\{a_0,a_2\}$ determine completely the sequence $\{a_j\}_{j\ge 0}$.
Also it is clear that when $j=0$ or~$2$, that matrix has nullity~$1$.
Therefore we can conclude that $\dim(S_\lambda)=2$.
The theorem follows.
\end{proof}

\begin{Theorem}
\label{pol}
Let~$H$ be the $\mathbb{C}^2$-valued analytic function on $(0,1)$ given by an irreducible sphe\-ri\-cal function~$\Phi$ on $
G$ of fundamental $K$-type $(1,\dots,1,0,\dots,0)\in\mathbb{C}^\ell$, with~$p$ ones, $0<p<\ell$.
If $P=\Psi^{-1}H$, then~$P$ is polynomial.
\end{Theorem}
\begin{proof}
We know that the function~$H$ is analytic in $(0,1)$, and from Corollary~\ref{operator2l} we know that it is an
eigenfunction of the operator $\widetilde D$ (see~\eqref{Dtilde}).
Also we know that the function $H(\tfrac{1+\cos s}{2})$ is analytic at $s=0$, since $\Phi(a(s))$ it is.
Therefore from Theorem~\ref{isomo} the function $P=\Psi^{-1}H$ is an analytic eigenfunction of~$D$ on the closed interval $[0,1]$.

If we introduce the following matrix-weight function $V=V(y)$ supported on the interval~$[0,1]$
\begin{gather*}
V(y)=\frac{(n-1)!!}{(n-2)!!}
\frac{2}{\omega_*} (y(1- y))^{n/2-1} \left(
\begin{matrix}
d_1&0
\\
0&d_2
\end{matrix}
\right),
\end{gather*}
with $\omega_*=\pi$ if~$n$ is even and $\omega_*=2$ if~$n$ is odd, then from Proposition~\ref{prodint} we have
\begin{gather*}
\langle \Phi_0,\Phi_1\rangle=\int_0^1H_2^*(y) V(y) H_1^*(y) dy.
\end{gather*}
It follows from Propositions~\ref{prodint} and~\ref{Deltasim} that $\widetilde D$ is a~symmetric operator
with respect to the inner product def\/ined among continuous vector-valued functions on $[0,1]$~by
\begin{gather*}
\langle H_1,H_2 \rangle_V =\int_0^1 H_2^*(y)V(y)H_1(y)dy.
\end{gather*}
Then, since $D=\Psi^{-1}\widetilde D \Psi$, we have that~$D$ is a~symmetric operator with respect to the inner product
def\/ined among continuous vector-valued functions on $[0,1]$~by
\begin{gather*}
\langle P_1,P_2 \rangle_W =\int_0^1 P_2^*(y)W(y)P_1(y)dy,
\end{gather*}
where
\begin{gather*}
W=\Psi^{*} V \Psi.
\end{gather*}
Actually, we have that $(W,D)$ is a~classical pair in the sense of~\cite{GPT03}, see also~\cite{D97}.
As the weight~$W$ has f\/inite moments there exists a~sequence $\{Q_r\}_{r\geq0}$ of $2\times 2$ matrix-valued orthonormal
polynomials, such that $DQ_r=Q_r\Lambda_r$ where $\Lambda_r$ is a~real diagonal matrix (for precise def\/initions and
general facts on matrix-valued orthogonal polynomials see~\cite{GPT02a} and~\cite{D97}).

Let $\{e_1,e_2\}$ be the canonical basis of $\mathbb{C}^2$.
Then
\begin{gather*}
\langle Q_re_j,Q_se_i\rangle_W=e_i^*\left(\int_0^1 Q_s^*(y)W(y)Q_r^*(y)dy \right)e_j=e_i^*\delta_{si}
Ie_j=\delta_{r,s}\delta_{i,j}.
\end{gather*}
Therefore, for $r\geq0$, $j=1, 2$, $\{Q_re_j\}$ is a~family of $\mathbb{C}^2$-valued orthonormal polynomials such that
\begin{gather*}
D(Q_re_j)=(DQ_r)e_j=(Q_r\Lambda_r) e_j=Q_r(\Lambda_r e_j)=\lambda_r^j(Q_r e_j),
\end{gather*}
where $\Lambda_r=\operatorname{diag}(\lambda_r^1,\lambda_r^2)$.

Now we write our function $P=\Psi^{-1}H$ as $P=\sum\limits_{r,j}a_{r,j}Q_re_j$, where $a_{r,j}=\langle P,Q_re_j\rangle_W$.
Since~$P$ is analytic on $[0,1]$ the sum converges not only in the $L^2$-norm but also in the topology based on uniform
convergence of sequences of functions and their successive derivatives.

Therefore,
\begin{gather*}
\lambda P=DP=\sum\limits_{r,j}a_{r,j}\lambda_r^jQ_re_j.
\end{gather*}
Then $a_{r,j}=0$ if $\lambda^j_{r}\neq\lambda$.
Since $\dim W_\lambda =2$ it follows that~$P$ is a~polynomial.
\end{proof}
\begin{Remark}
\label{remW}
It is easy to see from~\eqref{dimpar} and~\eqref{dimimpar} that the dimensions of the~$M$-submodules of the fundamental
representation of~$K$ with highest weight of the form $(1,\dots,1,0\dots,0)$, with~$p$ ones, are given~by
\begin{gather*}
d_1=\frac{(n-1)!}{(p-1)!(n-p)!},
\qquad
d_2=\frac{(n-1)!}{p!(n-1-p)!},
\end{gather*}
therefore the weight~$W$ is given~by
\begin{gather*}
W =\frac{(n-1)!!}{(n-2)!!}
\frac{2}{\omega_*} \frac{(n-1)!}{p!(n-p)!} (y(1- y))^{n/2-1} \Psi^* \left(
\begin{matrix}
p&0
\\
0&n-p
\end{matrix}
\right) \Psi,
\end{gather*}
with $\omega_*=\pi$ if~$n$ is even and $\omega_*=2$ if~$n$ is odd.
Then,~$W$ is a~scalar multiple of
\begin{gather*}
\left(
\begin{matrix}
p(2y-1)^2+n-p&n(2y-1)
\\
n(2y-1)&(n-p)(2y-1)^2+p
\end{matrix}
\right).
\end{gather*}
Even more, since $0<p<\ell$ and $n=2\ell,2\ell+1$ it follows that $p\neq n-p$.
Then it can be proved that the weight~$W$ does not reduce to a~smaller size, i.e., there is not any invertible
matrix~$M$ such that $M^*W(y)M$ is diagonal for all $y\in[0,1]$.
\end{Remark}

For a~given fundamental~$K$ type $\pi\in\hat{\mathrm{SO}}(n)$, $n=2\ell$ or $2\ell+1$, with highest weight of the form
$(1,\dots,1,0,\dots,0)\in\mathbb{C}^\ell$ with~$p$ ones ($0<p<\ell$), let $\Phi_{w,\delta}$ denote the irreducible
spherical function of the pair $({\mathrm{SO}}(n+1), {\mathrm{SO}}(n))$ given by $\tau\in\hat {\mathrm{SO}}(n+1)$ with
highest weight of the form
$(w+1,1,\dots,1,\delta,0\dots,0)$
with $p-1$ ones.

Therefore, combining~\eqref{lambda}, Theorems~\ref{polsol2l} and~\ref{pol} we have the following statement.

\begin{Theorem}
\label{columns}
Given $w\in\mathbb{N}_0$, every irreducible spherical function $\Phi_{w,\delta}$ of the pair
\mbox{$({\mathrm{SO}}(n+1),$} ${\mathrm{SO}}(n))$, with $n=2\ell$ or $2\ell+1$, of type
${\mathbf{m}}_n=(1,\dots,1,0,\dots,0)\in\mathbb{C}^\ell$ with~$p$ ones $(0<p<\ell)$, corresponds to a~vector valued
function $P_{w,\delta}$ $(\delta=0,1)$, which is a~polynomial of degree~$w$; and the leading coefficients of $P_{w,0}$
and $P_{w,1}$ are multiples of $\left(\begin{smallmatrix}
1
\\
0
\end{smallmatrix}\right)$ and $\left(\begin{smallmatrix}
0
\\
1
\end{smallmatrix}\right)$ respectively.
Precisely
\begin{gather*}
P_{w,\delta}(y)=\sum\limits_{j=0}^{w}\frac{y^j}{j!} [C;U;V+\lambda]_j P_{w,\delta}(0),
\end{gather*}
with
\begin{gather*}
C=
\begin{pmatrix}
(n/2+1)&1
\\
1&(n/2+1)
\end{pmatrix},
\qquad
U=(n+2)I,
\qquad
V= \left(
\begin{matrix}
p&0
\\
0&n-p
\end{matrix}
\right),
\\
\lambda=\lambda_{n}(w,\delta)=
\begin{cases}
-w(w+n+1)-p
&\text{if}
\quad
\delta=0,
\\
-w(w+n+1)-n+p
&\text{if}
\quad
\delta=1.
\end{cases}
\end{gather*}
Even more, the value of $P_{w,\delta}(0)$ can be computed.
\end{Theorem}
\begin{proof}
It only remains to prove that $P_{w,\delta}(0)$ can be computed.

Let us consider the case $\delta=0$.
We know from~\eqref{lambda} and Theorem~\ref{polsol2l} that there is some $c\in\mathbb{C}$ such that
\begin{gather*}
[C;U;V+\lambda]_w P_{w,0}(0)=c\left(\begin{matrix}
1
\\
0
\end{matrix}
\right).
\end{gather*}
Since $[C;U;V+\lambda]_w$ is invertible, this~$c$ is univocally determined by the condition $\Phi(e)=I$, which implies
\begin{gather*}
\Psi(1)\sum\limits_{j=0}^{w}\frac{1}{j!} [C;U;V+\lambda]_j P_{w,0}(0)=\left(\begin{matrix}
1
\\
1
\end{matrix}
\right).
\end{gather*}

Similarly, we can prove the same for $P_{w,1}(0)$.
\end{proof}

\begin{Remark}
\label{Ps}
It is worth to observe that for $w$, $w'\ge0$ and $\delta,\delta'=0,1$, since $\langle P_{w,\delta},P_{w',\delta'}
\rangle_W=\langle \Phi_{w,\delta},\Phi_{w',\delta'}\rangle$, we have that if $(w,\delta)\neq(w',\delta')$ then
\begin{gather*}
\langle P_{w,\delta},P_{w',\delta'} \rangle_W=0.
\end{gather*}
Therefore, our construction encodes all equivalent classes of irreducible spherical functions of a~fundamental $K$-type
of highest weight $\lambda_p$, $0<p<\ell$, in the orthogonal set of $\mathbb{C}^2$-valued polynomials
$\{P_{w,0},P_{w,1}\}$.
The degree of $P_{w,0}$ and $P_{w,1}$ is~$w$, and the leading coef\/f\/icient is a~multiple of $\left(\begin{smallmatrix}
1
\\
0
\end{smallmatrix}\right)$ or $\left(\begin{smallmatrix}
0
\\
1
\end{smallmatrix}\right)$, respectively.

\end{Remark}

\section{Matrix valued orthogonal polynomials}
\label{mvop}

\subsection{Matrix valued orthogonal polynomials}

In this subsection, given~$n$ of the form $2\ell$ or $2\ell+1$ with $\ell\in\mathbb{N}$, for a~f\/ixed $0< p<\ell$ we
shall construct a~sequence of matrix-valued polynomials $\{P_w\}_{w\ge0}$ directly related to irreducible spherical
functions of type $\pi\in\hat {\mathrm{SO}}(n)$ of highest weight
${\mathbf{m}}_\pi=(1,\dots,1,0\dots,0)\in\mathbb{C}^\ell$, with~$p$ ones.

Given a~nonnegative integer~$w$ and $\delta =0,1$, we can consider $\Phi_{w,\delta}$, the irreducible spherical function
of type~$\pi$ associated with the irreducible representation $\tau\in\hat {\mathrm{SO}}(n+1)$ of highest weight of the
form ${\mathbf{m}}_{\tau}=(w+1,1,\dots,1,\delta,0,\dots,0)$ with $p-1$ ones.

We insist on recalling that, since~$\pi$ has only two ${\mathrm{SO}}(n-1)$-submodules, we can interpret the diagonal
matrix-valued function $\Phi_{w,\delta}(a(s))$, $s\in(0,\pi)$, as a~$2$ column vector function.

Now we consider the vector-valued function
\begin{gather*}
P_{w,\delta}:\ (0,1)\to \mathbb{C}^2
\end{gather*}
given by the vector function $P_{w,\delta}(y)=\Psi^{-1}(y)\Phi_{w,\delta}(a(s))$, with $\cos(s)=2y-1$.
Then, we def\/ine the matrix-valued function
\begin{gather*}
P_w= P_w(y),
\end{gather*}
whose~$\delta$-th column ($\delta =0,1$) is given by the $\mathbb{C}^{2}$-valued polynomial $P_{w,\delta}(y)$.

Let consider the matrix-valued skew symmetric bilinear form def\/ined among $C^\infty$ $2\times 2$ matrix-valued functions
on $[0,1]$~by
\begin{gather*}
\langle P,Q \rangle_W =\int_0^1 Q^*(y)W(y)P(y)dy,
\end{gather*}
where
\begin{gather*}
W=\frac{(n-1)!!}{(n-2)!!}\frac{2}{\omega_*} \frac{(n-1)!}{p!(n-p)!} (y(1-y))^{n/2-1} \left(
\begin{matrix}
p(2y-1)^2+n-p&n(2y-1)
\\
n(2y-1)&(n-p)(2y-1)^2+p
\end{matrix}
\right).
\end{gather*}
See Remark~\ref{remW}.
Then we state the following theorem.
\begin{Theorem}
\label{MVOP}
The matrix-valued polynomial functions $P_w$, $ w\geq0$, form a~sequence of orthogonal polynomials with respect to $W$,
which are eigenfunctions of the symmetric differential operator~$D$ in~\eqref{D}.
Moreover,
\begin{gather*}
D P_w = P_w \left(\begin{matrix}
\lambda(w,0) & 0
\\
0 & \lambda(w,1)
\end{matrix}
\right),
\end{gather*}
where
\begin{gather*}
\lambda(w,\delta)=
\begin{cases}
-w(w+n+1)-p
&\text{if}
\quad
\delta=0,
\\
-w(w+n+1)-n+p
&\text{if}
\quad
\delta=1.
\end{cases}
\end{gather*}
\end{Theorem}

\begin{proof}
From Theorem~\ref{ecpar} we have that the~$\delta$-th column of $P_w$ is an eigenfunction of the operator~$D$ with
eigenvalue $\lambda(w,\delta)$, see~\eqref{lambda} and~\eqref{D}.
Therefore we have
\begin{gather*}
DP_w=P_w \Lambda_w,
\end{gather*}
with
\begin{gather*}
\Lambda_w=\left(
\begin{matrix}
\lambda(w,0) & 0
\\
0& \lambda(w,1)
\end{matrix}
\right).
\end{gather*}

From Theorem~\ref{columns} we know that each column of $P_w$ is a~polynomial function of degree~$w$ and, even more,
that $P_w$ is a~polynomial whose leading coef\/f\/icient is a~nonsingular diagonal matrix.

Given~$w$ and $w'$, non negative integers, by using Remark~\ref{Ps} we have
\begin{gather*}
\langle P_{w'},P_w \rangle_W = \int_{0}^1 P_w(y)^*W(y) P_{w'}(y)  du =\sum\limits_{\delta,\delta'=0}^1 \int_{0}^1
\big(P_{w,\delta}(y)^*W(y) P_{w',\delta'}(y) du \big) E_{\delta,\delta'}
\\
\phantom{\langle P_{w'},P_w \rangle_W}
=\sum\limits_{\delta,\delta'=0}^1 \delta_{w,w'}\delta_{\delta,\delta'} \left(\int_{0}^1 P_{w,\delta}(y)^*W(y)
P_{w',\delta'}(y) du \right) E_{\delta,\delta'}
\\
\phantom{\langle P_{w'},P_w \rangle_W}
=\delta_{w,w'}\sum\limits_{\delta=0}^1 \int_{0}^1 \big(P_{w,\delta}(y)^*W(y) P_{w'\delta}(y) du,\big)E_{\delta,\delta},
\end{gather*}
which proves the orthogonality.
Even more, it also shows us that $\langle P_w, P_{w} \rangle_W$ is a~diagonal matrix.
Also, making a~few simple computations we have that
\begin{gather*}
\langle D P_w,P_{w'}\rangle =\delta_{w,w'}\langle P_w,P_{w'}\rangle \Lambda_w =\delta_{w,w'} \Lambda_w^* \langle
P_w,P_{w'}\rangle =\langle P_w, DP_{w'}\rangle,
\end{gather*}
for every $w$, $w'\in\mathbb{N}_0$, since $\Lambda_w$ is real and diagonal.
This concludes the proof of the theorem.
\end{proof}

\section[The ${\mathrm{SO}}(2\ell+1)$-type with highest weight $2\lambda_\ell$]
{The $\boldsymbol{{\mathrm{SO}}(2\ell+1)}$-type with highest weight $\boldsymbol{2\lambda_\ell}$}
\label{PC}

In this section $K={\mathrm{SO}}(2\ell+1)$.
We will focus on the particular case when the $K$-type is given by an irreducible representation~$\pi$ with highest
weight $2\lambda_\ell=(1,1,\dots,1)$.
We will f\/irst see that such $K$-module is the direct sum of three~$M$-submodules, and we will f\/ind similar results to
those obtained for the fundamental $K$-types $\lambda_1,\dots,\lambda_{\ell-1}$ that are direct sum of
two~$M$-submodules.

Let us consider the irreducible $K$-module $\Lambda^\ell(V)$, with $V=\mathbb{C}^n$, $n=2\ell+1$.
The vector
$v=({\mathbf{e}}_1-i{\mathbf{e}}_2)\wedge({\mathbf{e}}_3-i{\mathbf{e}}_4)\wedge\dots\wedge({\mathbf{e}}_{2\ell-1}-i{\mathbf{e}}_{2\ell})$
is the unique, up to a~scalar, dominant vector and its weight is $2\lambda_\ell=(1,1,\dots,1)$.

It is not dif\/f\/icult to see that $\Lambda^\ell(V)$ is the sum of three~$M$-irreducible submodules, namely
\begin{gather}
\label{3M-mod}
\Lambda^\ell(V)=V_1\oplus V_0 \oplus V_{-1}
\end{gather}
with respective highest weights $(1,\dots,1), (1,\dots,1,0), (1,\dots,1,-1)\in\mathbb{C}^{\ell}$ and having
$V_0=\Lambda^{\ell-1}(V) \wedge{\mathbf{e}}_n$ and $V_1\oplus V_{-1}\simeq \Lambda^\ell(\mathbb{C}^{n-1})$.

The vectors
\begin{gather*}
v_1= ({\mathbf{e}}_1-i{\mathbf{e}}_2)\wedge({\mathbf{e}}_3-i{\mathbf{e}}_4)\wedge\dots\wedge({\mathbf{e}}_{2\ell-1}-i{\mathbf{e}}_{2\ell}),
\\
v_0= -({\mathbf{e}}_1-i{\mathbf{e}}_2)\wedge({\mathbf{e}}_3-i{\mathbf{e}}_4)\wedge\dots\wedge({\mathbf{e}}_{2\ell-3}-i{\mathbf{e}}_{2\ell-2})\wedge
{\mathbf{e}}_n,
\\
v_{-1}= ({\mathbf{e}}_1-i{\mathbf{e}}_2)\wedge({\mathbf{e}}_3-i{\mathbf{e}}_4)\wedge\dots\wedge({\mathbf{e}}_{2\ell-1}+i{\mathbf{e}}_{2\ell})
\end{gather*}
are~$M$-highest weight vectors in $ V_1$, $V_0$ and $ V_{-1}$, respectively.
Also let us call $P_1$, $P_0$ and $P_{-1}$ the respective projections on $ V_1$, $V_0$ and $ V_{-1}$, according to the
decomposition~\eqref{3M-mod}.

In order to obtain the explicit expression of~$E$ in~\eqref{vectoreq} we are interested to compute
\begin{gather*}
\sum\limits_{j=1}^{n-1}\dot\pi(I_{nj})P_{s}\dot\pi(I_{nj}){\big|_{V_{r}}}=\lambda(r,s)I_{V_{r}},
\end{gather*}
with $r,s=1,0,-1$ corresponding to the three~$M$-submodules $ V_1$, $V_0$ and $ V_{-1}$ of the representation~$\pi$.

If $1\le j\le\ell$, then
\begin{gather*}
\dot\pi(I_{n,2j-1})({\mathbf{e}}_{2k-1}-i{\mathbf{e}}_{2k})=
\begin{cases}
0&
\text{if}\quad
k\ne j,
\\
-{\mathbf{e}}_n &
\text{if}\quad
k=j,
\end{cases}
\\
\dot\pi(I_{n,2j})({\mathbf{e}}_{2k-1}-i{\mathbf{e}}_{2k})=
\begin{cases}
0&
\text{if}
\quad
k\ne j,
\\
i{\mathbf{e}}_n &
\text{if}\quad
k=j,
\end{cases}
\end{gather*}
therefore, it is easy to see that $P_0 \dot\pi(I_{n,2j-1}) v_0=P_0 \dot\pi(I_{n,2j}) v_0=0$ and that $P_r
\dot\pi(I_{n,2j-1}) v_s=P_r \dot\pi(I_{n,2j}) v_s=0$ when $s\pm1$ and $r\pm1 $; i.e.~
\begin{gather*}
\lambda(0,0)=\lambda(-1,-1)=\lambda(1,-1)=\lambda(-1,1)=\lambda(1,1)=0.
\end{gather*}
Furthermore, it is easy to see that, for $1\le j\le \ell$ and~$r$ equal to $1$ or $-1$, we have
\begin{gather*}
\dot\pi(I_{n,2j-1})P_0\dot\pi(I_{n,2j-1})v_r+\dot\pi(I_{n,2j})P_0\dot\pi(I_{n,2j})v_r=-v_r,
\end{gather*}
then $\lambda(-1,0)=\lambda(1,0)=-\ell$.
Therefore, it only remains to compute
\begin{gather*}
\sum\limits_{j=1}^\ell\left(\dot\pi(I_{n,2j-1})P_s\dot\pi(I_{n,2j-1})v_0+\dot\pi(I_{n,2j})P_s\dot\pi(I_{n,2j})v_0\right),
\end{gather*}
for $s=\pm1$.

To obtain $P_s\dot\pi(I_{n,k})v_0$ it is necessary to decompose $\dot\pi(I_{n,k})v_0$ according to the direct
sum~\eqref{3M-mod}.
We know that $\dot\pi(X_{-\varepsilon_j-\varepsilon_\ell})v_1\in V_1$ and
$\dot\pi(X_{-\varepsilon_j+\varepsilon_\ell})v_{-1}\in V_{-1}$; recall that
\begin{gather*}
X_{-\epsilon_j-\epsilon_\ell} =I_{2\ell-1,2j-1}-I_{2\ell,2j}+i(I_{2\ell-1,2j}+I_{2\ell,2j-1}),
\\
X_{-\epsilon_j+\epsilon_\ell} =I_{2\ell-1,2j-1}+I_{2\ell,2j}+i(I_{2\ell-1,2j}-I_{2\ell,2j-1}),
\end{gather*}
see~\eqref{rootvectors}.
We have
\begin{gather*}
\dot\pi\left(X_{-\varepsilon_j-\varepsilon_\ell}\right)\left({\mathbf{e}}_{2j-1}-i{\mathbf{e}}_{2j}\right)
=-2({\mathbf{e}}_{2\ell-1}+i{\mathbf{e}}_{2\ell}),
\\
\dot\pi\left(X_{-\varepsilon_j-\varepsilon_\ell}\right)\left({\mathbf{e}}_{2\ell-1}-i{\mathbf{e}}_{2\ell}\right)
=2\left({\mathbf{e}}_{2j-1}+i{\mathbf{e}}_{2j}\right),
\\
\dot\pi\left(X_{-\varepsilon_j-\varepsilon_\ell}\right)\left({\mathbf{e}}_{2k-1}-i{\mathbf{e}}_{2k}\right) =0,
\qquad
\text{for} \quad k\ne s,\ell.
\end{gather*}
Therefore, for $1\le j< \ell$,
\begin{gather*}
\dot\pi\left(X_{-\varepsilon_j-\varepsilon_\ell}\right)v_1=2({\mathbf{e}}_1-i{\mathbf{e}}_2)\wedge\dots
\wedge({\mathbf{e}}_{2(\ell-1)-1}-i{\mathbf{e}}_{2(\ell-1)})\wedge({\mathbf{e}}_{2j-1}+i{\mathbf{e}}_{2j})
\\
\qquad
{} -2({\mathbf{e}}_1-i{\mathbf{e}}_2)\wedge\dots\wedge{({\mathbf{e}}_{2j-3}-i{\mathbf{e}}_{2j-2})}
\wedge({\mathbf{e}}_{2\ell-1}+i{\mathbf{e}}_{2\ell})\wedge({\mathbf{e}}_{2j+1}-i{\mathbf{e}}_{2j+2})\wedge
\\
\qquad
\dots\wedge({\mathbf{e}}_{2\ell-1}-i{\mathbf{e}}_{2\ell}).
\end{gather*}
Similarly, for $1\le j< \ell$,
\begin{gather*}
\begin{split}
& \dot\pi\left(X_{-\varepsilon_j+\varepsilon_\ell}\right)v_1=2({\mathbf{e}}_1-i{\mathbf{e}}_2)
\wedge\dots\wedge({\mathbf{e}}_{2(\ell-1)-1}-i{\mathbf{e}}_{2(\ell-1)})\wedge({\mathbf{e}}_{2j-1}+i{\mathbf{e}}_{2j})
\\
& \qquad
{}+2({\mathbf{e}}_1-i{\mathbf{e}}_2)\wedge\dots\wedge{({\mathbf{e}}_{2j-3}-i{\mathbf{e}}_{2j-2})}
\wedge({\mathbf{e}}_{2\ell-1}+i{\mathbf{e}}_{2\ell})\wedge({\mathbf{e}}_{2j+1}-i{\mathbf{e}}_{2j+2})\wedge
\\
& \qquad
\dots\wedge({\mathbf{e}}_{2\ell-1}-i{\mathbf{e}}_{2\ell}).
\end{split}
\end{gather*}
Hence, for $1\le j< \ell$, we have
\begin{gather*}
\tfrac{-i}{8}\left(\dot\pi\left(X_{-\varepsilon_j-\varepsilon_\ell}\right)v_1+
\dot\pi\left(X_{-\varepsilon_j+\varepsilon_\ell}\right)v_{-1}\right)
\\
\qquad
=({\mathbf{e}}_1-i{\mathbf{e}}_2)\wedge\dots\wedge({\mathbf{e}}_{2(\ell-1)-1}-i{\mathbf{e}}_{2(\ell-1)})\wedge{\mathbf{e}}_{2j}
=\dot\pi (I_{n,2j})v_0,
\\
\tfrac{1}{8}\left(\dot\pi\left(X_{-\varepsilon_j-\varepsilon_\ell}\right)v_1+
\dot\pi\left(X_{-\varepsilon_j+\varepsilon_\ell}\right)v_{-1}\right)
\\
\qquad
=({\mathbf{e}}_1-i{\mathbf{e}}_2)\wedge\dots\wedge({\mathbf{e}}_{2(\ell-1)-1}-i{\mathbf{e}}_{2(\ell-1)})\wedge{\mathbf{e}}_{2j-1}
=\dot\pi (I_{n,2j-1})v_0.
\end{gather*}
Then, for $1\le j< \ell$,
\begin{gather*}
\dot\pi \left(I_{n,2j-1}\right)P_1\dot\pi \left(I_{n,2j-1}\right)v_0=\tfrac{1}{8}\dot\pi
\left(I_{n,2j-1}\right)\dot\pi\left(X_{-\varepsilon_j-\varepsilon_\ell}\right)v_1
\\
\qquad
=\tfrac
i2\big({\mathbf{e}}_1-i{\mathbf{e}}_2)\wedge\dots\wedge{\mathbf{e}}_{2j}\wedge
\dots\wedge({\mathbf{e}}_{2(\ell-1)-1}-i{\mathbf{e}}_{2(\ell-1)})\wedge{\mathbf{e}}_{n},
\\
\dot\pi \left(I_{n,2j}\right)P_1\dot\pi \left(I_{n,2j}\right)v_0=\tfrac{-i}{8}\dot\pi
\left(I_{n,2j}\right)\dot\pi\left(X_{-\varepsilon_j-\varepsilon_\ell}\right)v_1
\\
\qquad
=-\tfrac
12\big({\mathbf{e}}_1-i{\mathbf{e}}_2)\wedge\dots\wedge{\mathbf{e}}_{2j-1}\wedge
\dots\wedge({\mathbf{e}}_{2(\ell-1)-1}-i{\mathbf{e}}_{2(\ell-1)})\wedge{\mathbf{e}}_{n}.
\end{gather*}
Therefore, for $1\le j< \ell$,
\begin{gather*}
\dot\pi \left(I_{n,2j-1}\right)P_1\dot\pi \left(I_{n,2j-1}\right)v_0+\dot\pi \left(I_{n,2j}\right)v_0P_1\dot\pi
\left(I_{n,2j}\right)v_0=-\tfrac12v_0.
\end{gather*}
Besides, for $j=\ell$ we have
\begin{gather*}
\dot\pi (I_{n,2\ell})v_0=\tfrac{1}{2i}(-v_1+v_{-1})
\qquad
\text{and}
\qquad
\dot\pi (I_{n,2\ell-1})v_0=\tfrac{1}{2}(v_1+v_{-1}).
\end{gather*}
Therefore, since
\begin{gather*}
\dot\pi \left(I_{n,2\ell}\right)P_1\dot\pi (I_{n,2\ell})v_0 =-\tfrac{1}{2i}\dot\pi
\left(I_{n,2\ell}\right)v_1=-\tfrac12v_0,
\\
\dot\pi \left(I_{n,2\ell-1}\right)P_1\dot\pi (I_{n,2\ell-1})v_0 =\tfrac{1}{2}\dot\pi
\left(I_{n,2\ell-1}\right)v_1=-\tfrac12v_0,
\end{gather*}
we have that
\begin{gather*}
\sum\limits_{j=0}^{n-1} \dot\pi \left(I_{n,j}\right)P_1\dot\pi (I_{n,j})v_0=-\frac{\ell+1}{2}v_0,
\end{gather*}
i.e.
\begin{gather*}
\lambda(0,1)=-\frac{\ell+1}{2}.
\end{gather*}
Analogously we obtain
\begin{gather*}
\lambda(0,-1)=-\frac{\ell+1}{2}.
\end{gather*}
Hence
\begin{gather*}
(\lambda(r,s))_{-1\le r,s\le1}=
\begin{pmatrix}
0&-\ell&0
\\
-\tfrac{\ell+1}{2}&0&-\tfrac{\ell+1}{2}
\\
0&-\ell&0
\end{pmatrix}.
\end{gather*}

Therefore, we obtain a~more explicit version of Corollary~\ref{eigenfunction2} using~\eqref{vectoreq} and
Remark~\ref{ene}.
Confront Corollary~\ref{operator2l+1}.
\begin{Corollary}
\label{operator2l+12}
Let~$\Phi$ be an irreducible spherical function on~$G$ of type $\pi\in \hat {\mathrm{SO}}(n)$, $n=2\ell+1$.
If the highest weight of~$\pi$ is of the form $(1,\dots,1)\in\mathbb{C}^\ell$, then the function
$H:(0,1)\to\mathbb{C}^3$ associated with~$\Phi$ satisfies $\widetilde DH=\lambda H$, for some $\lambda\in\mathbb{C}$
with
\begin{gather*}
\widetilde DH=y(1-y)H''(y)+\frac{1}{2}n(1-2y)H'(y)+\frac{(1-2y)^2+1}{4y(1-y)}\left(\begin{matrix}
-\ell&0&0
\\
0&-\ell-1&0
\\
0&0&-\ell
\end{matrix}
\right)H(y)
\\
\phantom{\widetilde DH=}
{}+\frac{(1-2y)}{2y(1-y)}\left(\begin{matrix}
0&-\ell&0
\\
-\tfrac{\ell+1}{2}&0&-\tfrac{\ell+1}{2}
\\
0&-\ell&0
\end{matrix}
\right)H(y).
\end{gather*}
\end{Corollary}

\subsection[Spherical functions of ${\mathrm{SO}}(2\ell+1)$-type $2\lambda_\ell$]
{Spherical functions of $\boldsymbol{{\mathrm{SO}}(2\ell+1)}$-type $\boldsymbol{2\lambda_\ell}$}

Let $n=2\ell+1$, we now focus on the spherical functions $\Phi_{w,\delta}$ of type
${\mathbf{m}}_{n}=(1,\dots,1)\in\mathbb{C}^\ell$, which are associated with the irreducible representations of
${\mathrm{SO}}(n+1)$ of highest weights of the form ${\mathbf{m}}_{n+1}=(w+1,1,\dots,1,\delta) \in\mathbb{C}^{\ell+1}$
such that the following pattern holds
\begin{gather*}
\begin{matrix} w+1 &{} &1 &\dots &1 &{} &\delta &{}
\\
{} &1 &\dots &{} &\dots &1 &{} &-1
\end{matrix}.
\end{gather*}
As before we make the function~$\Psi$ whose columns are given by the spherical functions $\Phi_{0,\delta}$, $\delta=-1,0,1$.
When $w=0$, this is calculable using~\cite[p.~364, equation~(8)]{V92}
or alternatively by considering the~$G$-modules $\Lambda^{\ell+1}(\mathbb{C}^{n+1})=V_1\oplus V_{-1}$ and
$\Lambda^{\ell}(\mathbb{C}^{n+1})=V_0$ and working in the same way that we already did in the beginning of
Section~\ref{s2l} for the $2\times2$ cases (here $V_{t}$, for $t=1,0,-1$, are the irreducible~$G$-modules with highest
weights $(1,\dots,1,t)\in\mathbb{C}^{\ell+1}$).

Therefore, if $\cos s =2y-1$ we have
\begin{gather*}
\Psi(y) =\left(
\begin{matrix}
e^{i s}&1&e^{-i s}
\\
1& \tfrac12(e^{i s}+e^{-i s})&1
\\
e^{-i s}&1&e^{i s}
\end{matrix}
\right)\\
\hphantom{\Psi(y)}{}
 =\left(
\begin{matrix}
2y-1+2i\sqrt{y-y^2}&1&2y-1-2i\sqrt{y-y^2}
\\
1&2y-1&1
\\
2y-1-2i\sqrt{y-y^2}&1&2y-1+2i\sqrt{y-y^2}
\end{matrix}
\right).
\end{gather*}

Each column of~$\Psi$ satisf\/ies the dif\/ferential equation given in Corollary~\ref{operator2l+12}.
And it is easy to check that we have
\begin{gather*}
y(1-y)H''(y)+\frac{1}{2}n(1-2y)H'(y)+\frac{(1-2y)^2+1}{4y(1-y)}\left(\begin{matrix}
-\ell&0&0
\\
0&-\ell-1&0
\\
0&0&-\ell
\end{matrix}
\right)\Psi(y)
\\
\qquad
{}+\frac{(1-2y)}{2y(1-y)}\left(\begin{matrix}
0&-\ell&0
\\
-\tfrac{\ell+1}{2}&0&-\tfrac{\ell+1}{2}
\\
0&-\ell&0
\end{matrix}
\right)\Psi(y)=\Psi(y)\left(\begin{matrix}
-\ell-1&0&0
\\
0&-\ell&0
\\
0&0&-\ell-1
\end{matrix}
\right).
\end{gather*}

\begin{Theorem}
\label{ecimpar}
The function~$\Psi$ can be used to obtain a~hypergeometric differential equation from the one given in
Corollary~{\rm \ref{operator2l+12}}.
Precisely, if~$H$ is a~vector-valued solution of the differential equation in Corollary~{\rm \ref{operator2l+12}}, with
eigenvalue~$\lambda$, then $P=\Psi^{-1}H$ is a~solution of $DP=\lambda P$, where~$D$ is the hypergeometric differential
operator given~by
\begin{gather*}
DP=y(1-y)P''+(C-yU)P'-VP,
\end{gather*}
with
\begin{gather*}
C=\left(\begin{matrix}
(n+2)/2&1/2&0
\\
1&(n+2)/2&1
\\
0&1/2&(n+2)/2
\end{matrix}
\right),
\qquad
U=(n+2)I,
\\
V= \left(\begin{matrix}
-\ell-1&0&0
\\
0&-\ell&0
\\
0&0&-\ell-1
\end{matrix}
\right).
\end{gather*}
\end{Theorem}

\begin{proof}
Let us write $H=\Psi P$.
Then
\begin{gather*}
y(1-y)P''+\left(2y(1-y)\Psi^{-1}\Psi'+\frac n2(1-2y)I\right)P'
\\
\qquad
{}+\Psi^{-1}\left(y(1-y)\Psi''+\frac n2(1-2y)\Psi'+\frac{1+(1-2y)^2}{4y(1-y)}\left(\begin{matrix}
-\ell&0&0
\\
0&-\ell-1&0
\\
0&0&-\ell
\end{matrix}
\right)\Psi \vphantom{\left(\begin{matrix}
0&-\ell&0
\\
-\tfrac{\ell+1}{2}&0&-\tfrac{\ell+1}{2}
\\
0&-\ell&0
\end{matrix}
\right)} \right.
\\
\qquad
{}+\left.\frac{(1-2y)}{2y(1-y)}\left(\begin{matrix}
0&-\ell&0
\\
-\tfrac{\ell+1}{2}&0&-\tfrac{\ell+1}{2}
\\
0&-\ell&0
\end{matrix}
\right)\Psi\right)P=\lambda P.
\end{gather*}
Now we compute
\begin{gather*}
2y(1-y)\Psi^{-1}\Psi'+\frac n2(1-2y)I=-(n+2)yI+\left(\begin{matrix}
(n+2)/2&1/2&0
\\
1&(n+2)/2&1
\\
0&1/2&(n+2)/2
\end{matrix}
\right).
\end{gather*}
Therefore
\begin{gather*}
y(1-y)P''+\left(-(n+2)yI+\left(\begin{matrix}
(n+2)/2&1/2&0
\\
1&(n+2)/2&1
\\
0&1/2&(n+2)/2
\end{matrix}
\right)\right)P'
\\
\qquad
{}+\left(\left(\begin{matrix}
-\ell-1&0&0
\\
0&-\ell&0
\\
0&0&-\ell-1
\end{matrix}
\right)-\lambda I\right)P=0.
\end{gather*}
This completes the proof of the theorem.
\end{proof}

We obtain a~similar result to Theorem~\ref{polsol2l}, with an analogous proof:

\begin{Theorem}
\label{polsol2l+1}
For a~given $\ell\in\mathbb{N}$ let $n=2\ell+1$, then the nonzero polynomial eigenfunctions of
\begin{gather*}
DP=y(1-y)P''+(C-yU)P'-VP,
\end{gather*}
with
\begin{gather*}
C=\left(\begin{matrix}
(n+2)/2&1/2&0
\\
1&(n+2)/2&1
\\
0&1/2&(n+2)/2
\end{matrix}
\right),
\qquad
U=(n+2)I,
\\
V= \left(\begin{matrix}
-\ell-1&0&0
\\
0&-\ell&0
\\
0&0&-\ell-1
\end{matrix}
\right),
\end{gather*}
have eigenvalues $-w(w+n+1)-\ell$ or $-w(w+n+1)-\ell-1$, with $w\in\mathbb{N}_0$.
In both cases the degree of the polynomial is~$w$ and the leading coefficient can be any multiple of $\left(\begin{smallmatrix}
0
\\
1
\\
0
\end{smallmatrix}\right)$ or any linear combination of $\left(\begin{smallmatrix}
1
\\
0
\\
0
\end{smallmatrix}\right)$ and $\left(\begin{smallmatrix}
0
\\
0
\\
1
\end{smallmatrix}\right)$, respectively.

\end{Theorem}

Let us consider $\widetilde D $, the dif\/ferential operator on $(0,1)$ introduced in Corollary~\ref{operator2l+12}:
\begin{gather*}
\widetilde DH=y(1-y)H''(y)+\frac{1}{2}n(1-2y)H'(y)
\\
\phantom{\widetilde DH=}
{}+\frac{(1-2y)^2+1}{4y(1-y)}\left(\begin{matrix}
-\ell&0&0
\\
0&-\ell-1&0
\\
0&0&-\ell
\end{matrix}
\right)H(y) +\frac{(1-2y)}{2y(1-y)}\left(
\begin{matrix}
0&-\ell&0
\\
-\tfrac{\ell+1}{2}&0&-\tfrac{\ell+1}{2}
\\
0&-\ell&0
\end{matrix}
\right)H(y).
\end{gather*}
Recall that the operator~$D$ that appears in Theorem~\ref{polsol2l+1} extends the dif\/ferential operator $D=\Psi
\widetilde D \Psi^{-1}$ to the whole real line.

We want to focus our attention on the following vector spaces of $\mathbb{C}^{3}$-valued analytic functions on $(0,1)$:
\begin{gather*}
S_\lambda= \big\{H=H(y):\widetilde D H=\lambda H,~H(\tfrac{\cos s+1}{2})~\text{analytic at}~s=0\big\},
\\
W_\lambda= \big\{P=P(y):DP=\lambda P,~\text{analytic on}~[0,1]\big\}.
\end{gather*}

From Theorem~\ref{ecimpar} we know that the correspondence $P\mapsto\Psi P$ is an injective linear map from~$W_\lambda$
into~$S_\lambda$.
In fact, $\Psi((\cos s +1)/2)$ is analytic as a~function of~$s$ and~$P$ is analytic at $y=1$, hence $H((\cos s
+1)/2)=(\Psi P)((\cos s+1)/2)$ is analytic at $s=0$.

Then, we have an analogous result to Theorem~\ref{isomo}, whose proof is quite similar and therefore we will omit it.

\begin{Theorem}
The linear map $P\mapsto\Psi P$ is an isomorphism from $W_\lambda$ onto $S_\lambda$.
\end{Theorem}

Now, we can easily make a~proof similar to that one of Theorem~\ref{pol} in order to obtain next theorem.

\begin{Theorem}
\label{pol3x3}
Let~$H$ be the $\mathbb{C}^3$-valued analytic function on $(0,1)$ given by an irreducible sphe\-ri\-cal function~$\Phi$ on $
{\mathrm{SO}}(2\ell+2)$ of fundamental ${\mathrm{SO}}(2\ell+1)$-type $(1,\dots,1)\in\mathbb{C}^\ell$.
If $P=\Psi^{-1}H$, then~$P$ is polynomial.
\end{Theorem}

For a~given fundamental $K$-type $\pi\in\hat{\mathrm{SO}}(n)$, $n=2\ell+1$, with highest weight
$(1,\dots,1)\in\mathbb{C}^\ell$, let $\Phi_{w,\delta}$ denote the irreducible spherical function of the pair
$({\mathrm{SO}}(n+1), {\mathrm{SO}}(n))$ given by $\tau\in\hat {\mathrm{SO}}(n+1)$ with highest weight of the form
$(w+1,1,\dots,1,\delta)\in\mathbb{C}^{\ell+1}$, $\delta=-1,0,1$.

Now, combining Theorems~\ref{polsol2l+1},~\ref{pol3x3} and the expression of the eigenvalue $\lambda_n(w,\delta)$ given
in~\eqref{lambda} we have the following statement.

\begin{Theorem}
Given $w\in\mathbb{N}$, every irreducible spherical function $\Phi_{w,\delta}$ of the pair
\mbox{$({\mathrm{SO}}(n+1),$} ${\mathrm{SO}}(n))$ with $n=2\ell+1$, of type ${\mathbf{m}}_n=(1,\dots,1)\in\mathbb{C}^\ell$,
corresponds to a~vector-valued function~$P_{w,\delta}$ $(w\ge0$, $\delta=-1,0,1)$, which is a~polynomial of degree~$w$.
The leading coefficients of~$P_{w,0}$ is a~multiple of $\left(\begin{smallmatrix}
0
\\
1
\\
0
\end{smallmatrix}\right)$ and the leading coefficients of $P_{w,-1}$ and $P_{w,1}$ are both linear combinations of $\left(\begin{smallmatrix}
1
\\
0
\\
0
\end{smallmatrix}\right)$ and $\left(\begin{smallmatrix}
0
\\
0
\\
1
\end{smallmatrix}\right)$.
Precisely
\begin{gather*}
P_{w,\delta}(y)=\sum\limits_{j=0}^{w}\frac{y^j}{j!} [C;U;V+\lambda]_j P_{w,\delta}(0),
\end{gather*}
with
\begin{gather*}
C=\left(\begin{matrix}
(n+2)/2&1/2&0
\\
1&(n+2)/2&1
\\
0&1/2&(n+2)/2
\end{matrix}\right),
\\
U=(n+2)I,
\qquad
V= \left(\begin{matrix}
-\ell-1&0&0
\\
0&-\ell&0
\\
0&0&-\ell-1
\end{matrix}\right),
\\
\lambda=\lambda_{n}(w,\delta)=
\begin{cases}
-w(w+n+1)-\ell
&\text{if}
\quad
\delta=0,
\\
-w(w+n+1)-\ell-1
&\text{if}
\quad
\delta=\pm1.
\end{cases}
\end{gather*}
Even more, the value of $P_{w,\delta}(0)$ can be computed.
\end{Theorem}
\begin{proof}
It only remains to prove that $P_{w,\delta}(0)$ can be computed.

Let us consider the case $\delta=0$.
We know from~\eqref{lambda} and Theorem~\ref{polsol2l+1} that there is some $c\in\mathbb{C}$ such that
\begin{gather*}
[C;U;V+\lambda]_w P_{w,0}(0)=c\left(\begin{matrix}
0
\\
1
\\
0
\end{matrix}
\right).
\end{gather*}
Since $[C;U;V+\lambda]_w$ is invertible this~$c$ is univocally determined by the condition $\Phi(e)=I$ which implies
\begin{gather*}
\Psi(1)\sum\limits_{j=0}^{w}\frac{1}{j!} [C;U;V+\lambda]_j P_{w,0}(0)=\left(\begin{matrix}
1
\\
1
\\
1
\end{matrix}
\right).
\end{gather*}

Now let us consider the cases $\delta=\pm1$.
We know from~\eqref{lambda} and Theorem~\ref{polsol2l+1} that
\begin{gather*}
[C;U;V+\lambda]_w P_{w,\delta}(0)\in\left\langle \left(\begin{matrix}
1
\\
0
\\
0
\end{matrix}
\right),\left(\begin{matrix}
0
\\
0
\\
1
\end{matrix}
\right)\right\rangle;
\end{gather*}
since $[C;U;V+\lambda]_w$ is invertible, this condition tells us that $P_{w,\delta}(0)$ belongs to a~plane which
contains the origin and does not depend on~$\delta$.

Besides, the condition $\Phi_{w,\delta}(e)=I$, for $\delta=\pm1$, tells us
\begin{gather*}
\left(
\begin{matrix}
1
\\
1
\\
1
\end{matrix}
\right) =
\begin{pmatrix}
1&1&1
\\
1&1&1
\\
1&1&1
\end{pmatrix}
\sum\limits_{j=0}^{w}\frac{1}{j!} [C;U;V+\lambda]_j P_{w,\delta}(0).
\end{gather*}
Then, $P_{w,\delta}(0)$ belongs to a~plane, parallel to the kernel of
\begin{gather*}
\begin{pmatrix}
1&1&1
\\
1&1&1
\\
1&1&1
\end{pmatrix}
\sum\limits_{j=0}^{w}\frac{1}{j!} [C;U;V+\lambda]_j,
\end{gather*}
which does not contain the origin and does not depend on~$\delta$.
Therefore we know that both $P_{w,1}(0)$ and $P_{w,-1}(0)$ are in the same straight line.

On the other hand, recall that we have
\begin{gather*}
\Phi_{w,\delta}(a(s))=\Psi\left(\frac{\cos s+1}{2}\right) P_{w,\delta}\left(\frac{\cos s+1}{2}\right),
\end{gather*}
where
\begin{gather*}
a(s)=
\begin{pmatrix}
I_{n-1}&0&0
\\
0&\cos s&\sin s
\\
0&-\sin s&\cos s
\end{pmatrix}
,
\end{gather*}
then
\begin{gather*}
\frac{d}{ds}{\bigg|_{s=0}}\Phi(a(s))=
\begin{pmatrix}
i&0&-i
\\
0&0&0
\\
-i&0&i
\end{pmatrix}
P_{w,\delta}(1).
\end{gather*}
From~\cite[p.~364, equation~(8)]{V92} we can easily compute $\tfrac{d}{ds}\Phi_{w,\delta}(a(s))$ at $s=0$, which is obtained by looking at the
action of $\dot\tau(I_{n+1,n})$ and considering the corresponding projection, see~\eqref{proj}; having then
\begin{gather*}
\delta\frac{i(w+\ell+1)}{1+\ell}
\begin{pmatrix}
-1
\\
0
\\
1
\end{pmatrix}
=
\begin{pmatrix}
i&0&-i
\\
0&0&0
\\
-i&0&i
\end{pmatrix}
\sum\limits_{j=0}^{w}\frac{1}{j!} [C;U;V+\lambda]_j P_{w,\delta}(0).
\end{gather*}

This last condition establishes that $P_{w,1}(0)$ and $P_{w,-1}(0)$ are in two dif\/ferent and parallel planes, and the
line mentioned above does not belong to any of them since each plane has to intersect it.
Therefore the values of $P_{w,1}(0)$ and $P_{w,-1}(0)$ are univocally determined.
\end{proof}

\subsection{Matrix-valued orthogonal polynomials of size~3}

In this subsection, given~$n$ of the form $2\ell+1$ with $\ell\in\mathbb{N}$, we shall construct a~sequence of
matrix-valued polynomials $\{P_w\}_{w\ge0}$ directly related to irreducible spherical functions of type $\pi\in\hat
{\mathrm{SO}}(n)$ of highest weight ${\mathbf{m}}_\pi=(1,\dots,1)\in\mathbb{C}^\ell$.

Given a~nonnegative integer~$w$ and $\delta =-1,0,1$, we can consider $\Phi_{w,\delta}$, the irreducible spherical
function of type~$\pi$ associated with the irreducible representation $\tau\in\hat {\mathrm{SO}}(n+1)$ of highest weight
of the form ${\mathbf{m}}_{\tau}=(w+1,1,\dots,1,\delta)$.

We insist on recalling that, since~$\pi$ has only three ${\mathrm{SO}}(2\ell)$-submodules, we can interpret the diagonal
matrix-valued function $\Phi_{w,\delta}(a(s))$, $s\in(0,\pi)$, as a~$3$ column vector function.

Now we consider the vector-valued function
\begin{gather*}
P_{w,\delta}: \ (0,1)\to \mathbb{C}^3
\end{gather*}
given by the vector function $P_{w,\delta}(y)=\Psi^{-1}(y)\Phi_{w,\delta}(a(s))$, with $\cos(s)=2y-1$.
Then, we def\/ine the matrix-valued function
\begin{gather*}
P_w= P_w(y),
\end{gather*}
whose~$\delta$-th column ($\delta =-1,0,1$) is given by the $\mathbb{C}^{3}$-valued polynomial $P_{w,\delta}(y)$.

Let consider the matrix-valued skew symmetric bilinear form def\/ined among continuous $3\times 3$ matrix-valued functions
on $[0,1]$~by
\begin{gather*}
\langle P,Q \rangle_W =\int_0^1 Q^*(y)W(y)P(y)dy,
\end{gather*}
where the $3\times 3$ weight-matrix~$W$ is given~by
\begin{gather*}
W(y)=\frac{(n-1)!!}{(n-2)!!} (y(1- y))^{n/2-1} \Psi^*(y)\left(
\begin{matrix}
d_{1}&0&0
\\
0&d_2&0
\\
0&0&d_3
\end{matrix}
\right) \Psi(y)
\end{gather*}
with
\begin{gather*}
d_1=d_3=\frac{(2\ell+1)!}{\ell!(\ell+2)!},
\qquad
d_2=\frac{(2\ell+1)!}{\ell!(\ell+1)!},
\end{gather*}
and
\begin{gather*}
\Psi(y) =\left(
\begin{matrix}
2y-1+2i\sqrt{y-y^2}&1&2y-1-2i\sqrt{y-y^2}
\\
1&2y-1&1
\\
2y-1-2i\sqrt{y-y^2}&1&2y-1+2i\sqrt{y-y^2} \,
\end{matrix}
\right)
\end{gather*}
Let us recall that, from Proposition~\ref{prodint}, we have
\begin{gather*}
\langle \Phi_{w,\delta},\Phi_{w',\delta'}\rangle = \int_0^1 P_{w,\delta}^*W(y)P_{w',\delta'} dy.
\end{gather*}

\begin{Remark}
Notice that~$W$ reduces to a~smaller size: if $M=\left(\begin{smallmatrix}
1&0&1
\\
0&\sqrt2&0
\\
-1&0&1
\end{smallmatrix}\right) $ we have
\begin{gather*}
M W(y)M^*=\frac{(n-1)!!}{(n-2)!!} (y(1- y))^{n/2-1} 4
\\
\qquad{}
\times
\left(
\begin{matrix}
2d_1(2y-1)^2+d_2 & d_1(2y-1) + d_2 (2y-1)/\sqrt2 & 0
\\
d_1(2y-1)\sqrt2 + d_2 (2y-1)/\sqrt2 & d_1+d_2(2y-1)^2/2 & 0
\\
0&0&d_18(y-y^2)
\end{matrix}
\right).
\end{gather*}
\end{Remark}

Then we state the following theorem.
\begin{Theorem}
The matrix-valued polynomial functions $P_w$, $ w\geq0$, form a~sequence of ortho\-gonal polynomials with respect to~$W$,
which are eigenfunctions of the symmetric differential operator~$D$ from Theorem~{\rm \ref{ecimpar}}.
Moreover,
\begin{gather*}
D P_w = P_w \left(\begin{matrix}
\lambda(w,-1) & 0 & 0
\\
0 & \lambda(w,0) & 0
\\
0 & 0 & \lambda(w,1)
\end{matrix}
\right),
\end{gather*}
where
\begin{gather*}
\lambda(w,\delta)=
\begin{cases}
-w(w+n+1)-p
&\text{if}\quad
\delta=0,
\\
-w(w+n+1)-n+p
&\text{if}\quad
\delta=\pm1.
\end{cases}
\end{gather*}
\end{Theorem}

\begin{proof}
The proof is completely analogous to the proof of Theorem~\ref{MVOP}
\end{proof}

\section*{Appendix}
\pdfbookmark[1]{Appendix}{app}


\begin{proof}[Proof of Proposition~\ref{derprincipal}.]
For $|t|$ suf\/f\/iciently small $A(s,t)$ is close to the identity of~$K$, i.e.~to the identity mat\-rix~$I_{n}$.
So we can consider the function
\begin{gather}
\label{logA}
X(s,t)=\log(A(s,t))=B(s,t)-\frac{B(s,t)^2}2+\frac{B(s,t)^3}3-\cdots,
\end{gather}
where $B(s,t)=A(s,t)-I_{n}$.
Then
\begin{gather*}
\pi(A(s,t))=\pi(\exp X(s,t))=\exp\dot\pi (X(s,t))=\displaystyle\sum\limits_{j\geq 0}\frac{\dot\pi(X(s,t))^j}{j!}.
\end{gather*}
Now we dif\/ferentiate with respect to~$t$ to obtain
\begin{gather}
\label{derpi}
\frac{\partial(\pi\circ A)}{\partial t}=\dot\pi\left(\frac{\partial X}{\partial t}\right) +\tfrac
1{2!}\dot\pi\left(\frac{\partial X}{\partial t}\right)\dot\pi(X) +\tfrac 1{2!}\dot\pi(X) \dot\pi\left(\frac{\partial
X}{\partial t}\right)
\nonumber
\\
\phantom{\frac{\partial(\pi\circ A)}{\partial t}=}
{}+ \tfrac 1{3!}\dot\pi\left(\frac{\partial X}{\partial t}\right) \dot\pi(X)^2 + \tfrac 1{3!}\dot\pi (X)
\dot\pi\left(\frac{\partial X}{\partial t}\right) \dot\pi\left(X\right) + \tfrac 1{3!}\dot\pi(X)^2
\dot\pi\left(\frac{\partial X}{\partial t}\right) +\cdots.
\end{gather}
Since $X(s,0)=0$, if we dif\/ferentiate~\eqref{derpi} with respect to~$t$ and evaluate at $(s,0)$ we obtain
\begin{gather*}
\frac{\partial^2 (\pi\circ A)}{\partial t^2}\Big|_{t=0} = \dot \pi\left(\frac{\partial^2 X}{\partial
t^2}\Big|_{t=0}\right) +\dot\pi\left(\frac{\partial X}{\partial t}\Big|_{t=0}\right)^2.
\end{gather*}
To compute $\frac{\partial X}{\partial t}\big|_{t=0}$ and $\frac{\partial^2 X}{\partial t^2}\big|_{t=0}$ we
dif\/ferentiate~\eqref{logA} and we get
\begin{gather*}
\frac{\partial X}{\partial t}= \frac{\partial B}{\partial t} -\tfrac 1{2} \!\left(\frac{\partial B}{\partial t}\right)B-
\tfrac 12 B \!\left(\frac{\partial B}{\partial t}\right)+ \tfrac 13\! \left(\frac{\partial B}{\partial t}\right)B^2 +\tfrac 13
B\!\left(\frac{\partial B}{\partial t}\right) B +\tfrac 13 B^2 \!\left(\frac{\partial B}{\partial t}\right)+\cdots.
\end{gather*}
Since $B(s,0)=0$ we have
\begin{gather*}
\frac{\partial X}{\partial t}\Big|_{t=0}= \frac{\partial B}{\partial t}\Big|_{t=0}=\frac{\partial A}{\partial
t}\Big|_{t=0}.
\end{gather*}
We also get
\begin{gather*}
\frac{\partial^2 X}{\partial t^2}\Big|_{t=0} = \frac{\partial^2 A}{\partial t^2}\Big|_{t=0}-\left(\frac{\partial
A}{\partial t}\Big|_{t=0}\right)^2.
\end{gather*}

Now we will f\/irst consider the case $A(s,t)=k(s,t)$.
A~direct computation yields to
\begin{gather*}
\frac{\partial k}{\partial t}=
\begin{pmatrix}
{\bf 0} & {\bf 0} & {\bf 0} & {\bf 0} & {\bf 0}
\\
{\bf 0} & \frac{-\sin s\sin t}{(1-\cos^2s\cos^2t)^{3/2}} & {\bf 0} & \frac{\sin^2s\cos t}{(1-\cos^2s\cos^2t)^{3/2}} & 0
\\
{\bf 0} & {\bf 0} & {\bf 0} & {\bf 0} & {\bf 0}
\\
{\bf 0} & \frac{-\sin^2 s\cos t}{(1-\cos^2s\cos^2t)^{3/2}} & {\bf 0} & \frac{-\sin s\sin t}{(1-\cos^2s\cos^2t)^{3/2}} & 0
\\
{\bf 0} & 0 & {\bf 0} & 0 & 0
\end{pmatrix}
,
\end{gather*}
in particular $\frac{\partial k}{\partial t}\big|_{t=0}=\frac1{\sin s}I_{n,j}$.
Dif\/ferentiating once more with respect to~$t$ and evaluating at $t=0$ we obtain $\frac{\partial^2 k}{\partial
t^2}\big|_{t=0}=-\frac{1}{\sin^2 s}(E_{jj}+E_{n,n})$.
Then we get
\begin{gather*}
\frac{\partial^2 A}{\partial t^2}\Big|_{t=0}-\left(\frac{\partial A}{\partial t}\Big|_{t=0}\right)^2=-\frac{1}{\sin^2
s}(E_{jj}+E_{n,n})-\frac{1}{\sin^2 s}I_{n,j}^2=0.
\end{gather*}

Similarly when $A(s,t)=h(s,t)$ we obtain
\begin{gather*}
\frac{\partial h}{\partial t}=
\begin{pmatrix}
{\bf 0} & {\bf 0} & {\bf 0} & {\bf 0} & {\bf 0}
\\
{\bf 0} & \frac{-\sin s\cos^2 s\cos t\sin t}{(1-\cos^2s\cos^2t)^{3/2}} & {\bf 0} & \frac{-\cos s\cos t\sin^2s}{(1-\cos^2s\cos^2t)^{3/2}} & 0
\\
{\bf 0} & {\bf 0} & {\bf 0} & {\bf 0} & {\bf 0}
\\
{\bf 0} & \frac{\cos s\cos t\sin^2s}{(1-\cos^2s\cos^2t)^{3/2}} & {\bf 0} & \frac{-\sin s\cos^2s\cos t\sin t}{(1-\cos^2s\cos^2t)^{3/2}} & 0
\\
{\bf 0} & 0 & {\bf 0} & 0 & 0
\end{pmatrix},
\end{gather*}
in particular $\frac{\partial h}{\partial t}\big|_{t=0}=-\frac{\cos s}{\sin s}I_{n,j}$.
Dif\/ferentiating once more with respect to~$t$ and evaluating at $t=0$ we obtain $\frac{\partial^2 h}{\partial
t^2}\big|_{t=0}=-\frac{\cos^2s}{\sin^2 s}(E_{jj}+E_{n,n})$.
Then we get
\begin{gather*}
\frac{\partial^2 A}{\partial t^2}\Big|_{t=0}-\left(\frac{\partial A}{\partial t}\Big|_{t=0}\right)^2
=-\frac{\cos^2s}{\sin^2 s}(E_{jj}+E_{n,n})-\frac{\cos^2 s}{\sin^2 s}I_{n,j}^2=0.
\end{gather*}
Proposition follows.
\end{proof}

\subsection*{Acknowledgements}
This paper was partially supported by CONICET, PIP 112-200801-01533 and SeCyT-UNC.

\pdfbookmark[1]{References}{ref}
\LastPageEnding

\end{document}